\documentclass[11pt]{amsart}
\usepackage{amsmath}
\usepackage{amssymb}
\usepackage{euscript}
\usepackage{pifont}
\usepackage{tikz}
\usepackage[all]{xy}
\numberwithin{equation}{section} \numberwithin{figure}{section}
\input cyracc.def

\usepackage{amssymb,amsmath}
\usepackage{pstricks}
\usepackage{eepic}

\newrgbcolor{blue}{0 0 .7}

\newtheorem{definition}[equation]{Definition}
\newtheorem{theorem}[equation]{Theorem}
\newtheorem{proposition}[equation]{Proposition}

\newtheorem{lemma}[equation]{Lemma}
\newtheorem{example}[equation]{Example}
\newtheorem{remark}[equation]{Remark}

\newcommand{\w}{\wedge}
\newcommand{\n}{\notag}
\newcommand{\noi}{\noindent}
\newcommand{\im}{\mbox{i}}
\newcommand{\hook}{\hookrightarrow}

\newcommand{\vsp}{\vspace}
\newcommand{\sq}{ \;  \square}

\newcommand{\C}{ \mathbb{C} }

\newcommand{\PP}{ \mathbb{P} }
\newcommand{\del}{ \delta \phi}
\newcommand{\delo}{ \delta \omega}
\newcommand{\delt}{ \delta \theta}
\newcommand{\dele}{ \delta \eta}
\newcommand{\x}{\textnormal{x}}

\newcommand{\Pf} {\noi \emph{Proof.}$\; \;$}

\newcommand{\beq}{ \begin{equation}}
\newcommand{\eeq}{ \end{equation} }

\begin{document}

\sloppy

\title{Projectively   deformable  Legendrian surfaces}
\author{Joe S. Wang}
\email{jswang12@gmail.com}
\keywords{Legendrian surface,   projective deformation,  del Pezzo
surface} \subjclass[2000]{53A20}
\begin{abstract}
Consider  an immersed Legendrian surface in the five dimensional
complex projective space equipped with the standard homogeneous
contact structure. We introduce a class of fourth order projective
Legendrian deformation called \emph{$\,\Psi$-deformation}, and give
a differential geometric characterization of  surfaces admitting
maximum three parameter family of such deformations. Two explicit
examples of maximally $\, \Psi$-deformable surfaces are constructed;
the first one is given by a Legendrian map from $\, \PP^2$ blown up
at three distinct collinear points, which is an embedding away from
the -2-curve and degenerates to a   point along the -2-curve. The
second one is  a Legendrian embedding of  the degree 6 del Pezzo
surface, $\, \PP^2$ blown up at three non-collinear points. In both
cases, the Legendrian map is given by a system of cubics through the
three points, which is a subsystem of the anti-canonical system.
\end{abstract}
\maketitle \tableofcontents

\section{Introduction}\label{sec1}
Let $\, Z^{2m+1}$ be a complex manifold of  odd dimension $\, 2m+1
\geq 3$. A \emph{contact structure} on $\, Z$ is by definition a
hyperplane field $\, \mathcal{H} \subset T Z$ such that it  is
locally defined by $\, \mathcal{H} = \langle \,\alpha \,
\rangle^{\perp}$ for a  1-form $\, \alpha$ that satisfies the
nondegeneracy condition $\, \alpha \w (d \alpha)^m \ne 0$. A
manifold with a contact structure is called a \emph{contact
manifold}. Typical examples of  contact manifolds are the
homogeneous adjoint varieties of simple Lie algebras including the
odd dimensional projective spaces $\PP^{2m+1}$, \cite{LM}.

There exist a   distinguished class of subvarieties in a contact
manifold. A \emph{Legendrian subvariety} $\, M^m \subset Z^{2m+1}$
is a  $\, \mathcal{H}$-horizontal subvariety of maximum dimension
$\, m$. Typical examples of Legendrian subvarieties are the
homogeneous sub-adjoint varieties, \cite{LM}. In contrast to the
case of real and  smooth category, where both contact structures and
Legendrian subvarieties are flexible  and admit  unobstructed local
deformation, there are relatively small set of known Legendrian
subvarieties in complex contact manifolds, even in the simplest  odd
dimensional projective spaces.

Bryant showed that every compact Riemann surface can be embedded in
$\, \PP^3$ as a Legendrian curve in relation to the twistorial
description of minimal surfaces in 4-sphere,  \cite{Br1}. Landsberg
and Manivel adopted the idea from \cite{Br1} and showed that a K3
surface blown up at certain twelve points can be embedded in $\,
\PP^5$ as a Legendrian surface via an explicit  birational
contactomorphism from the projective cotangent bundle $\, \PP(T^*
\PP^3)$ to $\,\PP^5$. Buczynski established  a general principle of
hyperplane sections for Legendrian subvarieties in   projective
spaces, \cite{Bu2}. Successive hyperplane sections of the known
examples then gave many new   smooth  Legendrian subvarieties.
Buczynski also showed that the algebraic completion of the special
linear group $\, \overline{SL_3\C} \subset \PP^{17}$ is a smooth,
Fano Legendrian subvariety with  Picard number 1, \cite{Bu1}.

The purpose of the present paper is to propose a differential
geometric study of Legendrian surfaces in the five dimensional
projective space $\, \PP^5$ from the perspective of  projective
Legendrian deformation. Let us  explain the  motive. In his general
investigation of the deformation of a submanifold  in a homogeneous
space during the period 1916 -1920, Cartan considered the following
problem of third order projective deformation of a surface in the
projective space $\, \PP^3$, which built upon the earlier work of
Fubini, \cite{Ca2} and the reference therein; \emph{ Let
\textnormal{x}$: \Sigma \hook \PP^3$ be a surface. Let
\textnormal{x}$': \Sigma \hook \PP^3$ be a deformation of
\textnormal{x}.
 \textnormal{x}$'$ is a third order deformation of   \textnormal{x}
when there exists an application map $\, g: \Sigma \to SL_4 \C$ such
that for each $\, p \in \Sigma$,
 \textnormal{x} and $\, g(p) \circ $ \textnormal{x}$'$ agree up to order three at $\, p$.
Which surfaces in $\, \PP^3$  admit a nontrivial(application map $\,
g$ is nonconstant) third order deformation?}

Cartan showed that a generic surface is rigid and does not admit
such deformations, but that there exist two special  sets of
surfaces with nondegenerate second fundamental form that admit
maximum three parameter family of third order deformations.

The present work was  inspired especially by the observation that a
quartic Kummer surface is one of the maximally deformable surfaces,
\cite{Fe}. This led us to consider the analogous problem for
Legendrian surfaces in $\, \PP^5$, which  may allow one to obtain
the detailed structure equations for a set of Legendrian surfaces
with  special properties. The integration of the structure equation
so obtained may suggest a method of construction for   new examples.
In particular, one may hope to find the Legendrian analogue of
Kummer surfaces.

\vsp{0.5pc}

\textbf{Main results.}

1. Let  x$: M \hook \PP^5$ be a Legendrian surface in $\, \PP^5$
viewed as a homogeneous space of the symplectic group
$\, \textnormal{Sp}_3 \C $.
Assuming  the second fundamental form of x
is nondegenerate, the moving frame method is employed to determine
the basic local invariants of a Legendrian surface as  a set of
three symmetric differentials $\,(\, \Phi, \, \Psi, \, \chi)$ of
degree $\,(3, \, 3, \, 2)$ and order $\, (2, \, 4, \, 5)$,
Proposition \ref{2structure} and  \eqref{2Psii}.

2. A class of fourth order deformation called
\emph{$\Psi$-deformation} is introduced. Given a Legendrian surface
x$: M \hook \PP^5$,
 a deformation x$': M \hook \PP^5$ is a  $\Psi$-deformation
 when the pair  $(\, \Phi', \, \Psi')$ for  x$'$ is isomorphic  to that of  x
 at each point of $\, M$ up to motion by $\, \textnormal{Sp}_3 \C $.
This choice of deformation  is justified later by the analogy with
the aforementioned Cartan's  result that there exist two special
sets of Legendrian  surfaces called \emph{$\,D_0$-surfaces} and
\emph{$\,D$-surfaces} that admit maximum three parameter family of
$\, \Psi$-deformations, Proposition \ref{41pro}.

3.  The structure equations for the maximally $\, \Psi$-deformable
surfaces are determined. It turns out that  the local moduli space
of $\,D_0$-surfaces depends on 1 arbitrary function of one variable,
whereas the local moduli space of $\,D$-surfaces is finite
dimensional, Theorem \ref{4main}.

4.  Two global examples of  $\,D_0$-surfaces are constructed
explicitly.
The first one is given by a Legendrian map from
$\, \PP^2$ blown up at three distinct collinear points,
which is an embedding away from
the -2-curve(the proper transform of the line through
the three points) and degenerates to a  point along the -2-curve,
Theorem \ref{thm51}.
The second one is given by  a Legendrian embedding of
the degree 6 del Pezzo surface, $\, \PP^2$ blown up at three
non-collinear points, Theorem \ref{thm52}.

\vsp{0.5pc}

The paper is organized as follows. In Section \ref{sec2}, the basic
local invariants of a Legendrian surface   are defined as  a set of
three symmetric differentials $\,(\, \Phi, \, \Psi, \, \chi)$. By
imposing natural geometric conditions in terms of these invariants,
we identify four distinguished  classes of  Legendrian surfaces,
and determine their structure equations, Section  \ref{sec21}  through
\ref{sec24}. In Section \ref{sec3}, as a preparatory step for the
analysis of   deformation with geometric constraints in Section
\ref{sec4}, we  determine the structure equations for the second
order Legendrian deformation, or equivalently the Legendrian
deformation preserving the second order cubic differential $\,\Phi$.
The analysis shows that there is no local obstruction for the
second order deformation of a Legendrian surface with nondegenerate
cubic $\, \Phi$.
In Section \ref{sec4}, we introduce the $\Psi$-deformation,
which is the main object of study in this paper.
It is the second order deformation  which also preserves
the fourth  order cubic differential $\, \Psi$.
The analysis shows that
a generic Legendrian surface does not admit any nontrivial
$\,\Psi$-deformations,
but that the two special  sets of Legendrian surfaces
called $D_0$-surfaces  and  $D$-surfaces
admit maximum three parameter family of $\, \Psi$-deformations.
Moreover, a subset of these maximally deformable surfaces admit  $\,
\Psi$-deformations that also preserve the fifth order quadratic
differential $\, \chi$, Theorem \ref{43thm}.
In Section \ref{sec5},
we choose and integrate  two simple examples of  structure equations for
$\, D_0$-surfaces. In Section   \ref{sec51}, the flat case is
examined,  where all the structure coefficients vanish. The
structure equation is   integrated, and
one gets a rational Legendrian
variety with  a  single, second order branch type  isolated
singularity. Its  smooth resolution is $\, \PP^2$ blown up at three
distinct collinear points, with the exceptional divisor being the
\,-2-curve. In Section   \ref{sec52}, the case called
\emph{tri-ruled}(Section \ref{sec24}) is examined, where each leaf
of the three $\, \Phi$-asymptotic foliations lies in a linear
Legendrian $\, \PP^2$. The structure equation is integrated, and
one gets a smooth Legendrian embedding of $\, \PP^2$ blown up at
three non-collinear  points. The embedding is given by a subsystem
of the anti-canonical system and each of the six -1-curves is mapped
to a line.

Throughout the paper,
we freely apply the methods and results
of exterior differential systems.
We refer the reader to \cite{BCG} for the standard reference on the subject.

\section{Legendrian surface}\label{sec2}
The method of moving frames is a process of  equivariant frame
adaptation for a submanifold  in a homogeneous space. The
algorithmic operation of successive normalizations reveals the basic
local invariants of the submanifold as the coefficients of the
structure equation and their derivatives. The method was developed
by E. Cartan, and Cartan himself applied it extensively to a variety
of  problems.

In this section, the method of moving frames is applied to immersed
Legendrian surfaces in the five dimensional complex projective space
$\,\PP^5$.
We establish the fundamental structure equation
which depends on the sixth order jet of the Legendrian immersion, and
identify the basic local invariants as the set of  three symmetric
differentials $\, \Phi, \, \Psi$, and $\, \chi$ of order $\, 2, \,
4$, and $\, 5$ respectively, Proposition \ref{2structure}. $\, \Phi$
and $\, \Psi$ are cubic, and $\, \chi$ is quadratic.

In order to understand the geometric implication  of these
invariants, we impose a set of conditions in terms of
$\, \Phi, \,\Psi$, and $\, \chi$, and give  an analysis for
the Legendrian
surfaces that satisfy these conditions. This  in turn leads to  four
classes of Legendrian surfaces with interesting geometric
properties, Section \ref{sec21},  \ref{sec22}, \ref{sec23}, and
\ref{sec24} respectively.

Let us give an outline of the analysis.

$\bullet \; \Phi$,  cubic differential of order 2:\; $\, \Phi$
represents the second fundamental form of the Legendrian surface.
Assuming $\, \Phi$ is nondegenerate, the base locus of $\, \Phi$
defines a   3-web called \emph{asymptotic web}. It is the lowest
order local invariant of a  Legendrian surface.

We give an analysis for surfaces with flat asymptotic web, Section
\ref{sec21}. The condition for the asymptotic web to be flat is
expressed by a single fourth order equation for  the Legendrian
immersion, \eqref{21Kvanishing}. A differential analysis shows that
this PDE becomes involutive after a partial prolongation with the
general solution depending on five arbitrary functions of 1
variable, Proposition \ref{21pro}. It turns out that all of  the
surfaces that  are of interest to us necessarily have flat
asymptotic web, e.g., surfaces admitting maximum family of
nontrivial  Legendrian  deformation,  Section \ref{sec4}.

The moving frame computation  associates to each asymptotic
foliation a unique Legendrian $\, \PP^2$-field that has second order
contact with the given foliation, \eqref{22ruled}. An asymptotic
foliation is called \emph{ruled} when the associated Legendrian
$\, \PP^2$-field is leafwise constant. We give an analysis  for surfaces
with ruled asymptotic foliations, Section \ref{sec24}.

$\bullet \; \Psi$,   cubic differential of order 4:\; The moving
frame computations show that there is no third order local
invariants for a Legendrian surface.
The pencil of cubics $\, (  \Psi; \, \Phi )$ based at $\, \Phi$
accounts  for roughly one-half of the fourth order   invariants of a
Legendrian surface.

We give an analysis for surfaces with vanishing $\, \Psi$,  called
\emph{$\, \Psi$-null} surfaces, which can be considered as the
Legendrian analogue of  quadrics in $\, \PP^3$,    Section
\ref{sec22}. The condition for   $\, \Psi$  to vanish is expressed
by a pair of  fourth order equations  for  the Legendrian immersion,
\eqref{22defi}. A differential analysis shows that the structure
equation  for  $\, \Psi$-null  surfaces closes up with the general
solution depending on one constant, Proposition \ref{22pro}.

More generally, we give an analysis for the class of  surfaces
called  \emph{isothermally asymptotic} surfaces, which is the case
when $\, \Psi$ is proportional to $\, \Phi$ and the pencil
$\, ( \Psi; \, \Phi )$ degenerates, Section \ref{sec23}.
It will be shown
that  this class of surfaces are examples of surfaces admitting
maximum three parameter family of $\, \Psi$-deformations,
Section \ref{sec4}.

$\bullet \; \chi$-quadratic differential of order 5:\; The geometry
of $\, \chi$ is examined in Section \ref{sec4}. It will be shown
that there exist Legendrian surfaces which admit maximum one
parameter family of deformations preserving the triple
$\, (\Phi, \, \Psi, \, \chi)$.

For a modern exposition of Cartan's equivalence method,
we refer to
\cite{Ga}\cite{IL}.

\subsection{Structure equation}\label{20}
Let $\, V = \C^{6}$ be the six dimensional complex vector space. Let
$\, \varpi$ be the standard symplectic 2-form on $\, V$. Let $\,
\PP^5 = \PP(V)$ be the  projectivization equipped with the induced
contact structure. The contact hyperplane field $\, \mathcal{H}$ on
$\, \PP^5$ is defined by \beq \mathcal{H}_{\mbox{x}} = [ (
\hat{\mbox{x}} \lrcorner \, \varpi )^{\perp} ], \; \; \mbox{for}\;
\; \mbox{x} \in \PP^5,\n \eeq where   $ \hat{\mbox{x}} \in  V$ is
any de-projectivization of  x. $( \hat{\mbox{x}}  \lrcorner \,
\varpi )^{\perp} \subset V$ is a codimension one subspace
 containing $ \hat{\mbox{x}}$,
and its projectivization $\, [ ( \hat{\mbox{x}} \lrcorner \, \varpi
)^{\perp} ] \subset \PP^5$ is a hyperplane at x. The symplectic
group Sp$_3 \C$ acts transitively on $\, \PP^5$  as a  group of
contact transformation.

$\mathcal{H}$ inherits a conformal class of nondegenerate symplectic
2-form determined by the restriction of  $\, \varpi$ on the quotient
space $\, ( \hat{\mbox{x}}  \lrcorner \, \varpi )^{\perp} / \langle
\,    \hat{\mbox{x}} \,  \rangle$. A two dimensional Lagrangian
subspace of $\, \mathcal{H}_{\x}$ is called \emph{Legendrian}. Let
$\, \Lambda \to \PP^5$ be the bundle of Legendrian 2-planes. Let $\,
Lag( V)$ be the set  of  three dimensional Lagrangian subspaces of
$\, V$. The symplectic group Sp$_3 \C$ acts transitively on both $\,
\Lambda$ and $\, Lag( V)$, and there exists the incidence double
fibration;

\begin{picture}(300,109)(-67,-21) \label{2double}
\put(150,72){$\mbox{Sp$_3(\C)$}$} \put(172,57){$\pi$}
\put(162,57){$\downarrow$} \put(161,40){$\Lambda = \mbox{Sp$_3(\C)$}
/ P$} \put(175,25){$\searrow$} \put(145,25){$\swarrow$}
\put(184,10){ $Lag( V)$ } \put(129,10){ $ \PP^5$ }
\put(187,27){$\pi_1$} \put(133,27){$\pi_0$} \put(99, -11){
\textnormal{Figure \ref{2double}.  Double fibration }}
\end{picture}

\noi The fiber of $\, \pi_0$ is isomorphic to $\, Lag(2, \, \C^4) $,
and the fiber of $\, \pi_1$ is $\, \PP^2$.

To fix the notation once and for all, let us define the projection
maps $\, \pi,  \, \pi_0$, and  $\, \pi_1 \, $ explicitly. Let  $(e,
\, f) = (e_0, \, e_1, \, e_2, \, f_0, \, f_1, \, f_2)$ denote   the
Sp$_3\C \subset $ SL$_6\C$\, frame of $\, V$ such that  the 2-vector
$\, \varpi_{\flat} = e_0 \w f_0 + e_1 \w f_1+ e_2 \w f_2$ is dual to
the symplectic form $\, \varpi$. Define
\begin{align}\label{2defipi}
\pi (e, \,f) &= ([e_0], \, [e_0 \w e_1 \w e_2]),   \\
\pi_0 ([e_0], \, [e_0 \w e_1 \w e_2]) &= [e_0],  \n \\
\pi_1 ([e_0], \, [e_0 \w e_1 \w e_2]) &= [e_0 \w e_1 \w e_2]. \n
\end{align}
\noi In this formulation, the stabilizer subgroup $\, P$ in Figure
\ref{2double} is of the form \beq \label{2group} P  = \{ \,
\begin{pmatrix}
A       & B \\
\cdot  & (A^{t})^{-1}
\end{pmatrix} \},
\eeq where $\, (A^{-1} B) ^t = A^{-1} B$, and \beq A  = \{ \,
\begin{pmatrix}
*       & * & *  \\
\cdot & * & *   \\
\cdot & * & *
\end{pmatrix} \}. \n
\eeq Here '$\cdot$' denotes 0  and '$*$' is arbitrary.

The Sp$_3\C$-frame $\, (e, \, f)$  satisfies the structure equation
\beq\label{2frame} d  (e, \, f) = (e, \, f) \, \phi \eeq for the
Maurer-Cartan form $\, \phi$  of  Sp$_3\C$. The components of $\,
\phi$ are  denoted by \beq \phi =
\begin{pmatrix}
\omega &  \eta \\
\theta   &  - \omega^t
\end{pmatrix},\n
\eeq where $\, \{ \, \omega, \, \theta, \, \eta \, \}$ are 3-by-3
matrix 1-forms such that $\, \theta^t=\theta, \, \eta^t=\eta$. $\,
\phi$ satisfies the structure equation \beq\label{2MC} d \phi+\phi
\w \phi=0. \eeq

Let $\, \x : M \hook \PP^5$ be a $\, \mathcal{H}$-horizontal,
immersed Legendrian surface. We employ the method of moving frames
to normalize the Sp$_3 \C$-frame along $\, \x$. Our argument is
local, and the action of certain finite permutation group that
occurs  in the course of  normalization shall be ignored. This does
not affect the analysis nor the result of moving frame computation
for our purpose. The process of equivariant reduction terminates at
the sixth order jet of the immersion $\, \x$.

\texttt{1-adapted frame}. By definition, there exists a unique lift
$\, \tilde{\x}: M \hook \Lambda$. Let $\, \tilde{\x}^*$Sp$_3 \C \to
M$ be the pulled back $\, P$-bundle. We continue to use $\, \phi$ to
denote the  pulled back Maurer-Cartan form  on $\,
\tilde{\x}^*$Sp$_3 \C$. From  \eqref{2defipi}, \eqref{2frame}, the
initial state of $\, \phi$ on  $\, \tilde{\x}^*$Sp$_3 \C$ takes the
form \beq \phi =
\begin{pmatrix}
\omega_{{0 0}}&\omega_{{0 1}}&\omega_{{0 2}}
&\eta_{{0 0}}&\eta_{{0 1}}&\eta_{{0 2}} \\
\omega^{{1}}& \omega_{{1 1}}& \omega_{1 2} &\eta_{{1 0}}
&\eta_{{1 1}}&\eta_{{1 2}}\\
\omega^{{2}}&\omega_{2 1} & \omega_{2 2}
&\eta_{{2 0}}&\eta_{{2 1}}&\eta_{{2 2}}\\
\cdot & \cdot  & \cdot
&-\omega_{{0 0}}&-\omega^{{1}}&-\omega^{{2}}\\
\cdot &\theta_{1 1} &\theta_{1 2}
&-\omega_{{0 1}}&-\omega_{{1 1}}&-\omega_{{2 1}}\\
\cdot &            \theta_{2 1} &\theta_{2 2} &-\omega_{{0 2}}
&-\omega_{{1 2}}&-\omega_{{2 2}}
\end{pmatrix}, \n
\eeq where $\, \theta_{i j}=\theta_{j i}$,   $ \, \eta_{i j}=\eta_{j
i}$, and we denote $\, \omega_{i0} = \omega^i, \, i = 1, \, 2$. For
any section $s: \, M \to \tilde{\x}^*$Sp$_3(\C)$, $\, \{ \,
s^*\omega^1, \, s^*\omega^2 \}$ is a local coframe  of $\, M$.

\texttt{2-adapted frame}. Differentiating $\, \theta_{1 0}=0, \,
\theta_{2 0}=0$, one gets \beq
\begin{pmatrix}
\theta_{1 1}&\theta_{1 2} \\
\theta_{2 1}&\theta_{2 2}
\end{pmatrix} \w
\begin{pmatrix}
\omega^1 \\
\omega^2
\end{pmatrix} =0. \n
\eeq By Cartan's lemma, there exist coefficients $\, t_{ijk};  \,
i,j,k=1,2$, fully symmetric in   indices such that \beq \theta_{i j}
= t_{ijk} \, \omega^k. \n \eeq The structure equation shows that the
cubic differential \beq\label{2defiPhi} \Phi = \theta_{ij}\omega^i
\omega^j =  t_{ijk}  \, \omega^i   \omega^j  \omega^k \eeq is well
defined on $\, M$  up to  scale. $\, \Phi$ represents   the second
fundamental form of  the Legendrian immersion.
\begin{definition}
Let $\, \x: M \hook \PP^5$ be an immersed Legendrian surface. Let
$\, \Phi$ be the  cubic differential  \eqref{2defiPhi} which
represents  the second fundamental form of the immersion $\, \x$.
The Legendrian surface  is  \emph{nondegenerate} if the cubic
differential  $\, \Phi$ is equivalent to an element in the unique
open orbit  of the general linear group GL$_2 \C$ action on cubic
polynomials in two variables, \cite{Mc}.
\end{definition}
\begin{remark}
The Segre embedding $\, \PP^1\times Q^1 \subset \PP^5$ is ruled by
lines, and it has a degenerate second fundamental cubic.
\end{remark}

We assume the Legendrian surface  is  nondegenerate from now on. By
a frame adaptation, one may   normalize $\, t_{ijk}$ such that
\begin{align}\label{2Phi}
\Phi &=  3  \, \omega^1   \omega^2 (\omega^1+\omega^2),  \\
       &= -3  \, \omega^1   \omega^2  \omega^3, \n
\end{align}
where $\, \omega^3=- (\omega^1 + \omega^2)$. This is equivalent to
\beq\label{2theta}
\begin{pmatrix}
\theta_{1 1}&\theta_{1 2} \\
\theta_{2 1}&\theta_{2 2}
\end{pmatrix} =
\begin{pmatrix}
\omega^2                    & \omega^1+\omega^2  \\
\omega^1+\omega^2  & \omega^1
\end{pmatrix}.
\eeq The structure group $\, P$, \eqref{2group}, for the 2-adapted
frame is  reduced  such that \beq A  = \{ \,
\begin{pmatrix}
a       & *   \\
\cdot & a\, A'
\end{pmatrix} \; | \quad   a \ne 0  \; \}, \n
\eeq where $\, A'$ is the finite  subgroup of  GL$_2 \C$ whose
induced action  leaves $\, \Phi$ invariant.

Three asymptotic line fields are determined  by $\, \{ \,
(\omega^1)^{\perp}, \, (\omega^2)^{\perp}, \, (\omega^3)^{\perp} \,
\}$. The set of respective  foliations defines    a 3-web called
\emph{asymptotic web}  on the Legendrian surface. Since a planar
3-web has local invariants, e.g.,  web curvature, asymptotic web is
the lowest order  invariant of  a nondegenerate Legendrian surface.

\texttt{3-adapted frame}. On the 2-adapted frame satisfying
\eqref{2theta}, set \beq\label{2setomega}
 \omega_{ij} = \frac{\delta_{ij}}{3} \, \omega_{00} + s_{ijk} \, \omega^k,
\quad    \mbox{for} \; \;  i,j =1, 2, \eeq for coefficients $\,
s_{ijk}$. Differentiating \eqref{2theta}, one gets
\begin{align}\label{2sfirst}
-3\,s_{{212}}+s_{{221}}+2\,s_{{111}}+2\,s_{{211}}&=0,   \\
-\frac{1}{2}\,s_{{121}}+\frac{3}{2}\,s_{{112}}+\frac{1}{2}\,s_{{212}}
-\frac{3}{2}\,s_{{221}}+s_{{211}}-s_{{122}}&=0, \n \\
-3\,s_{{121}}+s_{{112}}+2\,s_{{122}}+2\,s_{{222}}&=0. \n
\end{align}
Exterior derivatives of  \eqref{2setomega} show that
\begin{align}
ds_{112} &\equiv - \eta_{11}-\eta_{21}+\frac{1}{3} \omega_{02}, \n \\
ds_{211} &\equiv -\eta_{22}, \n \\
ds_{212} &\equiv - \eta_{21}-\eta_{22}+                     \omega_{01}, \n \\
ds_{121} &\equiv - \eta_{11}-\eta_{21}+                     \omega_{02}, \n \\
ds_{122} &\equiv -\eta_{11}, \n \\
ds_{221} &\equiv - \eta_{21}-\eta_{22}+\frac{1}{3} \omega_{01},
\quad \mod   \; \omega^1, \, \omega^2; \, \omega_{00}.  \n
\end{align}
By  a  frame adaptation, one may  translate the coefficients
$\,\{ s_{211}=0, \, s_{122}=0, \, 3  s_{112}-s_{121}=0  \}$,
which forces
$\, 3 s_{221}-s_{212}=0$ by \eqref{2sfirst}.
This set of normalizations
is equivalent to adapting the Sp$_3\C$-frame so  that
\begin{align}
d e_1 &\equiv f_2 \, \omega^1 \mod \; e_0, \, e_1; \, \omega^2, \n \\
d e_2 &\equiv f_1 \, \omega^2 \mod \; e_0, \, e_2; \, \omega^1, \n \\
d (e_1-e_2) &\equiv (f_1+f_2) \, \omega^2 \mod \; e_0, \, (e_1-e_2);
\, \omega^3. \n
\end{align}
\noi By  a  further frame adaptation, one may  translate
$\, s_{ijk}=0$(we omit the details), and we have \beq\label{2omega}
 \omega_{ij} = \frac{\delta_{ij}}{3} \, \omega_{00}.
\eeq

For this  3-adapted frame,  the structure equation shows that a
triple of    Legendrian $\, \PP^2$-fields is well defined along the
Legendrian surface. Let $\, (L_1, \, L_2, \, L_3)$ be the   triple
of   Legendrian $\, \PP^2$-fields, or equivalently the  triple of
Lagrangian 3-plane fields, defined by
\begin{align}\label{2Lagtriple}
(L_1, \, L_2, \, L_3)&:  \;M \to Lag(V)\times Lag(V) \times Lag(V),   \\
L_1 &= [e_0 \w e_1 \w f_2],  \n  \\
L_2 &= [e_0 \w e_2 \w f_1], \n \\
L_3 &= [e_0 \w (e_1-e_2)  \w (f_1+f_2)]. \n
\end{align}
Each  $\, L_i$ is the unique Legendrian $\, \PP^2$-field that has second
order contact with the asymptotic foliation defined by $\, \langle
\omega^i   \rangle^{\perp}$.

\texttt{4-adapted frame}. On the 3-adapted frame satisfying
\eqref{2omega}, set
\begin{align}\label{2eta}
 \omega_{01} &= h_{01k} \, \omega^k,   \\
 \omega_{02} &= h_{02k} \, \omega^k, \n \\
 \eta_{ij} =\eta_{ji}        &= h_{ijk} \, \omega^k,
     \quad \mbox{for} \; \;  i,j =1, 2. \n
\end{align}
The structure equation shows that the cubic differential
$\, \Psi =\eta_{ij}  \omega^i   \omega^j$ is well defined up to scale,
and up to translation by
\beq \Psi   \to  \Psi
+ (s_1 \, \omega^1+s_2\,\omega^2)   ((\omega^1)^2+(\omega^2)^2),
\quad \mbox{for  arbitrary coefficients} \; s_1, \, s_2. \n
\eeq
By a frame adaptation, one may translate
$\, h_{111}=0, \, h_{222}=0$ so that
\beq\label{2Psi}
\Psi = \omega^1   \omega^2
  ( (h_{112}+2 \, h_{121})\, \omega^1 + (2 \,h_{122}+h_{221})\, \omega^2 ).
\eeq
\noi $\Psi$ is now  well defined up to scale. It  is a fourth
order invariant of the nondegenerate Legendrian surface. For the
problem of   projective deformation of Legendrian surfaces, $\,
\Psi$ will play the role of the third fundamental form for surfaces
in $\, \PP^3$.

Note that the derivative   of  \eqref{2omega} with the relation
$\, h_{111}=0, \, h_{222}=0$ gives  the compatibility equations
\begin{align}\label{2etarela}
h_{011}&=h_{121}+h_{221},   \\
h_{022}&=h_{112}+h_{122}, \n \\
h_{012}&= \frac{3}{5}\, ( - h_{{12 1}}+ h_{{12 2}} ), \n \\
h_{021}&= \frac{3}{5} \, (  h_{{12 1}} - h_{{12 2}} ). \n
\end{align}

\texttt{5-adapted frame}. On the 4-adapted frame satisfying
\eqref{2Psi}, \eqref{2etarela}, set
\begin{align}\label{2eta0}
 \eta_{10} &= h_{10k} \, \omega^k,  \\
 \eta_{20} &= h_{20k} \, \omega^k. \n
\end{align}
The structure equation shows that the quadratic differential
$\,\chi = \eta_{10} \, \omega^1 +\eta_{20} \, \omega^2$
is well defined up to scale, and up to translation by
\beq
\chi   \to  \chi + s_0 \,((\omega^1)^2+(\omega^2)^2),
\quad \mbox{for  arbitrary coefficient}
\; s_0. \n
\eeq
By a frame adaptation, one may translate $\,
h_{101}+ h_{202}=0$(we omit the details)  so that \beq\label{2chi}
\chi = h_{101} (\omega^1)^2+ (h_{102}+h_{201}) \omega^1  \omega^2
         -h_{101} (\omega^2)^2.
\eeq \noi $\chi$ is now  well defined up to scale. It  is a fifth
order invariant of the nondegenerate Legendrian surface.

Restricting to the sub-bundle defined by the equation  $\, h_{101}+
h_{202}=0$, we set
\beq\label{2eta00}
\eta_{00} = h_{001}\,\omega^1+ h_{002}\, \omega^2.
\eeq
At this stage,  no more frame
adaptation is available, and the components of
the induced Maurer-Cartan form $\, \phi$ is uniquely determined,
modulo at most a finite
group action (this finite group does not enter into our analysis,
and we shall not pursue the exact expression for the representation
of this  group). The reduction process of  moving frame method
stops here.

For a notational purpose, let us make a change of variables;
\begin{align}\label{2notation}
h_{112}&=a_1, \quad \quad h_{101}=b_1, \quad \quad h_{001}=c_1,   \\
h_{121}&=a_2, \quad \quad h_{102}=b_2, \quad \quad h_{002}=c_2.\n \\
h_{122}&=a_3, \quad \quad h_{201}=b_3,  \n \\
h_{221}&=a_4,   \n
\end{align}
 The  covariant derivatives are denoted by
\begin{align}
d a_i &= -\frac{4}{3} \, a_i \,\omega_{00}  + a_{ik}\, \omega^k, \n \\
d b_i &= - 2                 \, b_i \,\omega_{00} \;+ b_{ik}\, \omega^k, \n \\
d c_i &= -\frac{8}{3} \, c_i \,\omega_{00} + c_{ik}\, \omega^k. \n
\end{align}
Differentiating \eqref{2eta}, \eqref{2eta0}, \eqref{2eta00}, one
gets a set of compatibility equations among the covariant
derivatives.
\begin{align}\label{2relation}
a_{11}&=-2\, b_2,  \\
a_{{22}}&=-\frac{3}{2}\,a_{{2 1}}+\frac{15}{2} \,b_{{3}}+8\,b_{{1}}, \n  \\
a_{{31}}&=-\frac{3}{2} \,a_{{21}}+\frac{15}{2} \,b_{{3}}+10\,b_{{1}}, \n \\
a_{{32}}&=a_{{2 1}}+5\,b_{{2}}-12\,b_{{1}}-5\,b_{{3}}, \n \\
a_{{4 2}}&=-2 \,b_{{3}}, \n \\
b_{{1 2}}&=b_{{2
1}}-\frac{3}{5}\,{a_{{3}}}^{2}+c_{{2}}+a_{{4}}a_{{1}}
                     -\frac{2}{5}\,a_{{2}}a_{3}, \n \\
b_{{3 2}}&=-b_{{1 1 }}-c_{{1}}+\frac{3}{5}\,{a_{{2}}}^{2}
                   +\frac{2}{5}\,a_{{2}}a_{{3}}-a_{{4}}a_{{1}}, \n \\
c_{{1 2}}&=c_{{2
1}}-2\,a_{{1}}b_{{3}}-2\,a_{{3}}b_{{3}}+2\,a_{{2}}b_{{2}}
                       +2\,a_{{4}}b_{{2}}. \n
\end{align}

The structure equation \eqref{2MC} for $\, \phi$ is now an identity
with these relations. One may  check that the  structure equation
with  the   coefficients $\, \{ \, a_i, \, b_j, \, c_k\, \}$ becomes
involutive after one prolongation with the general solution
depending on one arbitrary function of 2 variables in the sense of
Cartan, \cite{BCG}.
\begin{proposition}\label{2structure}
Let $\, \x: M \hook \PP^5$ be a nondegenerate  immersed Legendrian
surface. Let $\, \tilde{\x}: M \hook \Lambda$ be the associated lift
to the bundle of Legendrian 2-planes. Let $\, \tilde{\x}^*$Sp$_3(\C)
\to M$ be the pulled back bundle,   Figure \ref{2double}. $\,
\tilde{\x}^*$Sp$_3(\C)$ admits a reduction to a sub-bundle with
1-dimensional fibers such that the induced Maurer-Cartan form $\,
\phi$ satisfies the structure equations \eqref{2theta},
\eqref{2omega}, \eqref{2eta}, \eqref{2Psi}, \eqref{2etarela},
\eqref{2eta0}, \eqref{2eta00}, \eqref{2notation}, and
\eqref{2relation}.
\end{proposition}
\noi We shall work with the 5-adapted frame for the rest of the
paper. Unless stated otherwise, 'the structure equation' would mean
the structure equation for  the 5-adapted frame.

Note that  under the notations we chose, the invariant differentials
$\, \Phi$, $\, \Psi$,  and $\, \chi$, \eqref{2Phi}, \eqref{2Psi},
\eqref{2chi}, are expressed by
\begin{align}\label{2Psii}
\Phi &= 3\, (\omega^1)^2  \omega^2 +  3\, \omega^1  ( \omega^2)^2,  \\
\Psi  &= \omega^1   \omega^2  ( (a_1+2 \, a_2)\, \omega^1
          + (2 \,a_3+a_4)\, \omega^2 ),  \n  \\
      &= (a_1+2 \, a_2)\, (\omega^1)^2   \omega^2
          + (2 \,a_3+a_4)\,\omega^1  ( \omega^2)^2, \n \\
\chi &=  b_1 (\omega^1)^2+ (b_2+b_3) \omega^1  \omega^2
             -b_1 (\omega^2)^2. \n
\end{align}
In the next two sections, Section \ref{sec3} and  Section
\ref{sec4}, we shall examine the deformability, or the rigidity, of
Legendrian surfaces preserving these invariant differentials.

Before we proceed to the problem of deformation, let us examine four
classes of Legendrian surfaces with special geometric  properties.
There exist a  number of  surfaces in $\, \PP^3$ with notable
characteristics, which have been the subject of  extensive study,
\cite{Fe}. Some of the surfaces described  below can be considered
as the Legendrian analogues of these classical surfaces.

\subsection{Flat asymptotic 3-web}\label{sec21}
In this sub-section, we consider the class of Legendrian surfaces
with flat asymptotic 3-web. For a comprehensive introduction to web
geometry, we refer to \cite{PP}.

\begin{definition}
Let $\, M \hook \PP^5$ be a nondegenerate Legendrian surface. The
\emph{asymptotic 3-web}     is the set of three   foliations
defined by $\, \{ \,  (\omega^1)^{\perp}, \, (\omega^2)^{\perp}, \,
(\omega^3)^{\perp} \, \}$ at  2-adapted frame, where $\,
\omega^1+\omega^2+\omega^3=0$.
\end{definition}
\noi The following analysis shows  that the differential equation
describing the Legendrian surfaces with flat asymptotic web is in
good form(involutive) and admits arbitrary function worth solutions
locally.

The web curvature of the asymptotic 3-web can be expressed in terms
of the structure coefficients of the Legendrian surface. From the
structure equation(for  5-adapted frame), \beq d \omega^i =
\frac{2}{3}\, \omega_{00} \w  \omega^i, \quad i = 1, \, 2, \, 3. \n
\eeq The web curvature $\, K$ of the 3-web  is given by
\begin{align}\label{21webK}
 \frac{2}{3}\, d \omega_{00} &=K \, \omega^1 \w \omega^2,   \\
                                                  &= - \frac{4}{5} \, (a_2 - a_3)  \omega^1 \w \omega^2.  \n
\end{align}
The asymptotic web is flat  when \beq\label{21Kvanishing} a_2 -
a_3=0. \eeq We wish to give an  analysis of the compatibility
equations derived from this vanishing condition.

Differentiating $\, a_2 - a_3=0$, one gets
\begin{align}
a_{21}&=-2\,b_1+3\,b_2, \n  \\
b_3&=-2\,b_1+b_2. \n
\end{align}
Differentiating the second equation for $\, b_3$, one gets
\begin{align}
b_{22}&=2\,b_{{2 1}}-b_{{1 1}}-{a_{{2}}}^{2}+a_{{4}}a_{{1}}-c_{{1}}+2\,c_{{2}}, \n \\
b_{31}&=-2\,b_{11}+b_{21}. \n
\end{align}
Exterior derivative $\, d(d(a_2)) =0$ with these relations  then
gives \beq b_{21}=- b_{{1
1}}+{a_{{2}}}^{2}-a_{{4}}a_{{1}}+3\,c_{{1}}-4\,c_{{2}}. \n \eeq The
identities from $\, d(d(b_1)) =0, \, d(d(b_2)) =0$ determine the
derivative of $\, b_{11}$ by
\begin{align}\label{21b11}
d b_{11}  +\frac{8}{3} \, b_{11} \, \omega_{00} &=
 (-a_{{4}}b_{{2}}-\frac{1}{4} \,a_{{4 1}}a_{{1}}-2\,b_{{1}}a_{{2}}
 -2\,b_{{1}}a_{{1}}+b_{{2}}a_{{1}}+\frac{1}{4} \,a_{{1 2}}a_{{4}}+c_{{2 2}}
 +2\,c_{{1 1}}-3\,c_{{2 1}}) \,\omega^1   \\
 &\,+(a_{{4}}b_{{2}}+\frac{1}{4} \,a_{{4 1}}a_{{1}}+2\,b_{{1}}a_{{2}}
 +2\,b_{{1}}a_{{1}}-b_{{2}}a_{{1}}-\frac{1}{4} \,a_{{1 2}}a_{{4}}-c_{{2 2}}+c_{{1 1}}
) \,\omega^2. \n
\end{align}

At this step, we interrupt the differential analysis and invoke a
version of   Cartan-K\"ahler theorem, a general existence theorem
for  analytic differential systems, \cite{BCG}.
\begin{proposition}\label{21pro}
The structure equation for the nondegenerate Legendrian surfaces
with flat asymptotic web is in involution with the general solution
depending on five arbitrary functions of 1 variable.
\end{proposition}
\Pf From the analysis above, the exterior derivative identities $\,
d(d(a_1)) =0$, $\, d(d(a_4)) =0$, $\, d(d(c_1)) =0$, $\, d(d(c_2))
=0$, $\, d(d(b_{11})) =0$ give 5 compatibility equations while the
remaining  independent derivative coefficients at this step are $\,
\{ \, a_{12},  \,a_{41}, \, c_{11}, \, c_{21}, \, c_{22}; \, b_{11}
\, \}$. An inspection shows that the resulting structure equation is
in involution with the last nonzero Cartan character s$_1=5$. $\sq$

\subsection{Vanishing cubic differential $\, \Psi$}\label{sec22}
In this sub-section, we consider the class of Legendrian surfaces
with  vanishing cubic differential $\, \Psi$, a fourth order
invariant   \eqref{2Psi}.
\begin{definition}
Let $\, M \hook \PP^5$ be a nondegenerate Legendrian surface. $\, M$
is  a \emph{$\, \Psi$-null  surface} if the fourth order cubic
differential $\, \Psi$ defined at 4-adapted frame vanishes.
\end{definition}
\noi The following analysis shows  that a $\, \Psi$-null surface
necessarily has flat asymptotic web, and that the local moduli space
of $\, \Psi$-null surfaces is finite dimensional.

From \eqref{2Psii}, $\, \Psi$ vanishes when \beq\label{22defi}
a_1=-2\,a_2, \; \; a_4=-2\,a_3. \eeq Differentiating these
equations, one gets
\begin{align}
a_{21}&=  b_2,          \hspace{1.7cm}     a_{12}=-15\,b_3-16\,b_1+3\,b_2, \n   \\
b_3&=-2\,b_1+b_2,    \quad      a_{41}= 10\,b_1 - 12\,b_2. \n
\end{align}
Exterior derivatives $\, d(d(a_2)) =0, \, d(d(a_3)) =0$ with these
relations give
\begin{align}
b_{22}&= 8\,b_{{1 1}}+\frac{8}{5}
\,{a_{{2}}}^{2}-\frac{8}{5}\,a_{{2}}a_{{3}}
+\frac{15}{2}\,b_{{3 1}}-\frac{3}{2}\,b_{{2 1}},  \n \\
b_{31}&=-\frac{6}{5} \,b_{{1 1}}+b_{{2 1}}-{\frac
{88}{125}}\,{a_{{3}}}^{2} +\frac{8}{5}\,c_{{2}}+{\frac
{196}{125}}\,a_{{2}}a_{{3}}- \frac{6}{5}\,c_{{1}} +{\frac
{42}{125}}\,{a_{{2}}}^{2}. \n
\end{align}

Differentiating $\, b_3 =-2\,b_1+b_2$, one gets
\begin{align}
b_{11}&={\frac {22}{25}}\,{a_{{3}}}^{2}-2\,c_{{2}}
-{\frac {49}{25}}\,a_{{2}}a_{{3}}
+\frac{3}{2}\,c_{{1}}-{\frac {21}{50}}\,{a_{{2}}}^{2}, \n \\
b_{21}&=2\,c_{{1}}-{\frac {22}{25}}\,{a_{{2}}}^{2}-{\frac
{41}{25}}\,a_{{2}}a_{{3}} +{\frac
{51}{50}}\,{a_{{3}}}^{2}-\frac{5}{2} \,c_{{2}}. \n
\end{align}
Exterior derivatives $\, d(d(b_1)) =0, \, d(d(b_2)) =0$ with these
relations give
\begin{align}
c_{{2 2}}&=c_{{1 1}}+{\frac {68}{5}}\,b_{{2}}a_{{3}}+{\frac
{58}{5}}\,b_{{1}}a_{{2}}
-{\frac {68}{5}}\,a_{{2}}b_{{2}}-{\frac {78}{5}}\,b_{{1}}a_{{3}}, \n \\
c_{{2 1}}&=\frac{4}{3} \,c_{{1 1}}+{\frac
{266}{15}}\,b_{{2}}a_{{3}}+6\,b_{{1}}a_{{2}}
-18\,b_{{1}}a_{{3}}-{\frac {32}{5}}\,a_{{2}}b_{{2}}. \n
\end{align}

The identities from $\, d(d(c_1)) =0, \, d(d(c_2)) =0$ determine the
derivative of $\, c_{11}$ by
\begin{align}
dc_{11}+\frac{10}{3} \,c_{{1 1}}\omega_{{0 0}}&=
  ( {\frac {6807}{875}}\,{a_{{2}}}^
{3}-{\frac {607}{35}}\,c_{{2}}a_{{3}}-6\,{b_{{2}}}^{2}+{\frac
{9993}{ 875}}\,{a_{{3}}}^{3}+20\,b_{{1}}b_{{2}}-{\frac
{6646}{875}}\,{a_{{2}}}
^{2}a_{{3}}-12\,{b_{{1}}}^{2} \n \\
&\quad +{\frac {638}{35}}\,c_{{1}}a_{{3}}-{\frac
{513}{35}}\,c_{{1}}a_{{2}} -{\frac {12779}{875}}\,a_{{2}}{a_{{3}}
}^{2}+{\frac {447}{35}}\,c_{{2}}a_{{2}} ) \, \omega^{{1}} \n \\
&\quad + ( {\frac {10119}{875}}\,{a_{{2}}}^{3}-{\frac
{856}{35}}\,c_{{2}}a_{{3}}+ 62\,{b_{{2}}}^{2}+{\frac
{13016}{875}}\,{a_{{3}}}^{3}-118\,b_{{1}}b_{{
2}}-{\frac {8747}{875}}\,{a_{{2}}}^{2}a_{{3}}+54\,{b_{{1}}}^{2} \n \\
&\quad +{\frac {165}{7}}\,c_{{1}}a_{{3}}-{\frac
{141}{7}}\,c_{{1}}a_{{2}} -{\frac
{11763}{875}}\,a_{{2}}{a_{{3}}}^{2}+{\frac
{771}{35}}\,c_{{2}}a_{{2}} )\, \omega^{{2}}. \n
\end{align}
Differentiating this equation again, one gets a compatibility
equation of the form \beq\label{22c11} (a_2-a_3)\, c_{11} =  [ \,
a_i, \, b_j, \, c_k \,], \eeq where the right hand side is  a
polynomial in the variables  $\, a_i, \, b_j, \, c_k $. At this
juncture, the analysis  divides into two cases.

\texttt{Case $\, a_2 - a_3\ne 0$}. It turns out that the condition
$\, a_2 - a_3\ne 0$ is not compatible with the vanishing of $\,
\Psi$, and there is no nondegenerate Legendrian surfaces with $\,
\Psi\equiv0$, and $\, a_2 - a_3\ne 0$. Some of the expressions for
the analysis of this case are long. Let us explain the relevant
steps of differential analysis, and omit the details of the long and
non-essential terms.

From \eqref{22c11},  solve for $\, c_{11}$. Differentiating this,
one gets a set of  two equations, from which one  solves for $\,
c_1, \, c_2$. Differentiating these equations, one gets another set
of two equations which imply $\, b_1=b_2=0$. Differentiating these
equations again, one finally gets two quadratic equations for $\,
a_2, \, a_3$, which force  $\, a_2=a_3=0$, a contradiction.

\texttt{Case $\, a_2 - a_3=0$}. From Section \ref{sec21}, this is
the case when the asymptotic 3-web is flat. Successive derivatives
of the equation $\, a_2-a_3=0$ imply the following.
\begin{align}
b_1&=b_2,  \n \\
c_1&=c_2,  \n \\
c_{11}&=2\,a_2  b_2. \n
\end{align}
Furthermore, these equations are compatible, i.e., $\, d^2=0$ is an
identity.

The  remaining independent coefficients at this step are  $\, \{ \,
a_2, \, b_2, \, c_2 \, \}$. Let us remove the sub-script, and denote
$\, \{\, a_2, \, b_2, \, c_2 \,\}= \{\, a, \, b, \, c \,\}$. The
structure equations for these coefficients  are reduced to
\begin{align}\label{22dabc}
d a&=  -\frac{4}{3} a\,  \omega_{00} +  b (\omega^1-\omega^2),   \\
d b&= -2\, b\, \omega_{00}
- \frac{1}{2} (3\,{a}^{2} + c  )  (\omega^1-\omega^2), \n \\
d c&= -\frac{8}{3} \,  c\, \omega_{00} +  2\,ab\,
(\omega^1-\omega^2). \n
\end{align}
The Maurer-Cartan form $\, \phi$  takes the form \beq\label{22Mphi}
\phi =
\begin{pmatrix}\omega_{{0 0}}&-a\omega^{{1}}&-
a \omega^{{2}} &\;\;\;c( \omega^{{1}} +\omega^{{2}}) &b
(\omega^{{1}} + \omega^{{2}})&-b ( \omega^{{1}} +\omega^{{2}})
\\\noalign{\medskip}\omega^{{1}}&\frac{1}{3}\,\omega_{{0 0}}
&\cdot&\;\;\;b(\omega^{{1}} +\omega^{{2}})&-2\,a \omega^{{2}}
&\;\;\;a (\omega^{{1}}+ \omega^{{2}}) \\
\noalign{\medskip}\omega^{{2}}&\cdot&\frac{1}{3}\,\omega_{{0 0}}
&-b(\omega^{{1}}+\omega^{{2}})
&a(\omega^{{1}}+ \omega^{{2}}) &-2\,a\omega^{{1}}\\
\noalign{\medskip}\cdot&\cdot&\cdot&-\omega_{{0 0}}&-\omega^{{1}}&-\omega^{{2}}\\
\noalign{\medskip}\cdot&\omega^{{2}}&\omega^{{1}}+\omega^{{2}}&a\omega^{{1}}
&-\frac{1}{3}\,\omega_{{0 0}}&\cdot \\
\noalign{\medskip}\cdot&\omega^{{1}}+\omega^{{2}}&\omega^{{1}}&a
\omega^{{2}} &\cdot&-\frac{1}{3}\,\omega_{{0 0}}
\end{pmatrix}.
\eeq
\begin{proposition}\label{22pro}
Let $\, M \hook \PP^5$ be a nondegenerate, $\, \Psi$-null Legendrian
surface. The asymptotic 3-web  of   $\, M$ is necessarily flat. The
Maurer-Cartan form of the 5-adapted frame of $\, M$ is reduced to
\eqref{22Mphi}, and the structure coefficients $\, \{ \, a, \, b, \,
c \, \}$ satisfy the   equation \eqref{22dabc}. The local moduli
space of  $\, \Psi$-null Legendrian surfaces  has general dimension
1.
\end{proposition}
\Pf Let  F$\to M$ be the canonical bundle of   5-adapted frames
from Proposition \ref{2structure}. \eqref{22dabc} shows that the
invariant map $\, (a, \, b, \, c ):$ F$\to \C^3$ generically has
rank two. From the general theory of geometric structures with
closed structure equation, \cite{Br2}, the local moduli space of
this class of Legendrian surfaces has general dimension  \
dim$(\C^3)-$rank$(a,b,c)=1$. A Legendrian surface in this class
necessarily possesses a  minimum  1-dimensional  local group of
symmetry. The line field $\, \langle \, \omega_{00}, \,
\omega^1-\omega^2 \rangle^{\perp}$ is  tangent to the fibers of the
invariant map  $\, (a, \, b, \, c )$, and it generates a   local
symmetry. $\sq$

\subsection{Isothermally asymptotic}\label{sec23}
In this sub-section, we consider the class of Legendrian surfaces
which are the    analogues of the classical isothermally asymptotic
surfaces in $\, \PP^3$, \cite{Fe}.
\begin{definition}
Let $\, M \hook \PP^5$ be a nondegenerate Legendrian surface. $\, M$
is \emph{isothermally asymptotic} if  the fourth order cubic
differential $\, \Psi$, \eqref{2Psii}, is a multiple of the second
order cubic differential $\, \Phi$, \eqref{2Phi}.
\end{definition}
\noi The following analysis shows  that the differential equation
describing the isothermally asymptotic Legendrian surfaces is in
good form and admits arbitrary function worth solutions  locally.

From  \eqref{2Psii}, $\, \Psi \equiv  0 \mod \Phi$ when \beq
a_1+2\,a_2 = 2\, a_3+a_4. \n \eeq We wish to give an  analysis of
the compatibility equations derived from this condition.

Differentiating $\, a_1+2\,a_2 = 2\, a_3+a_4$, one gets
\begin{align}
a_{12}&=5\,a_{{2 1}}+10\,b_{{2}}-27\,b_{{3}}-40\,b_{{1}}, \n  \\
a_{41}&=5\,a_{{2 1}}-2\,b_{{2}}-15\,b_{{3}}-20\,b_{{1}}. \n
\end{align}
The identities from $\, d(d(a_1)) =0, \, d(d(a_2)) =0$ determine the
derivative of $\, a_{21}$ by
\begin{align}
d a_{21}  +2\,a_{{2 1}} \, \omega_{00} &= ( -2\,b_{{2 1}}+8\,b_{{1
1}}-\frac{2}{5} \,b_{{2 2}}+{\frac {27}{5}}\,b_{{3 1}}
-{\frac {8}{25}}\,a_{{2}}a_{{1}}+{\frac {8}{25}}\,a_{{1}}a_{{3}} ) \omega^{{1}} \n \\
&\quad+ ( -4\,b_{{1 1}}+\frac{3}{5}\,b_{{2 2}}+3\,b_{{2 1}}
-\frac{3}{5}\,b_{{3 1}}+\frac{8}{5}\,{a_{{2}}}^{2}
-\frac{8}{5}\,a_{{2}}a_{{3}}+{\frac {12}{25}}\,a_{{2}}a_{{1}}
-{\frac {12}{25}}\,a_{{1}}a_{{3}} ) \omega^{{2}}. \n
\end{align}
Exterior derivative $\, d(d(a_3)) =0$ with these relations  then
gives \beq b_{22}=5\,b_{{1 1}}+b_{{3 1}}+5\,b_{{2 1}}+{\frac
{46}{5}}\,a_{{2}}a_{{1}} -{\frac
{46}{5}}\,a_{{1}}a_{{3}}-15\,c_{{1}}-{\frac {44}{5}}\,{a_{{3}}}^{2}
+20\,c_{{2}}+5\,{a_{{1}}}^{2}+{\frac
{21}{5}}\,{a_{{2}}}^{2}-\frac{2}{5}\,a_{{2}}a_{{3}}. \n \eeq

At this step, we interrupt the differential analysis and invoke a
version of  Cartan-K\"ahler theorem.
\begin{proposition}
The structure equation for the nondegenerate isothermally asymptotic
Legendrian surfaces is in involution with the general solution
depending on five arbitrary functions of 1 variable.
\end{proposition}
\Pf From the analysis above, the   exterior derivative identities
$\, d(d(b_1)) =0$, $\, d(d(b_2)) =0$, $\, d(d(b_3)) =0$, $\,
d(d(c_1)) =0$, $\, d(d(c_2)) =0$, $\, d(d(a_{21})) =0$ give 6
compatibility equations while the remaining  independent derivative
coefficients at this step are $\, \{ \, b_{11},  \,b_{21}, \,
b_{31}, \, c_{11}, \, c_{21}, \, c_{22}; \, a_{21}  \, \}$. A short
analysis shows that the resulting structure equation becomes
involutive after one prolongation with  the last nonzero Cartan
character s$_1=5$. Since the prolonged structure equation does not
enter into our analysis in later sections, the details shall be
omitted. $\sq$


\subsubsection{Isothermally asymptotic with flat asymptotic web}\label{sec231}
Consider the class of isothermally asymptotic Legendrian surfaces
which have    flat asymptotic 3-web. From \eqref{2Psii} and
\eqref{21webK}, this  is equivalent to the condition
\beq\label{231defi} a_1=a_4, \quad a_2=a_3. \eeq The following
analysis shows that the differential equation describing  such
Legendrian surfaces is still in good form and admits arbitrary
function worth solutions  locally. Note that a $\, \Psi$-null
surface is necessarily isothermally asymptotic with flat asymptotic
web.

This is in contrast with the $\, \Psi$-null surface case, where the
defining equation \eqref{22defi} is also a set of two linear
equations among   $\, a_i$'s and yet the resulting structure
equations close up to admit solutions with finite dimensional
moduli. This reflects the subtle well-posedness of the equation
\eqref{231defi}. The discovery of this class of Legendrian surfaces
is perhaps most unexpected of the analysis in this section.

We wish to give an analysis for  the compatibility equations derived
from  \eqref{231defi}.
Differentiating  the given equations $\, a_1=a_4, \, a_2=a_3$,
one gets
\begin{align}
a_{21} &=4\,b_1+3\,b_3, \quad  \; \,   b_3 =-2\,b_1+b_2, \n  \\
a_{12} &=4\,b_1-2\,b_2, \quad    a_{41} =-2\,b_2.\n
\end{align}
Exterior derivatives  $\, d(d(a_1)) =0, \, d(d(a_2)) =0, \,
d(d(a_3)) =0$  with these relations give
\begin{align}
b_{22}&= b_{31}, \n \\
b_{31}&=-\frac{5}{3} \,b_{{1 1}}+\frac{4}{3}\,b_{{2 1}}+\frac{1}{3}
\,{a_{{1}}}^{2}
-\frac{1}{3}\,{a_{{2}}}^{2}-c_{{1}}+\frac{4}{3}\,c_{{2}}, \n \\
b_{21}&=-b_{{1
1}}+3\,c_{{1}}-4\,c_{{2}}+{a_{{2}}}^{2}-{a_{{1}}}^{2}. \n
\end{align}

Successively differentiating  $\, b_3 =-2\,b_1+b_2$, one gets
\begin{align}
c_1&=c_2, \n \\
c_{21}&=c_{11}, \quad c_{22} = c_{{1
1}}+4\,b_{{1}}a_{{1}}+4\,b_{{1}}a_{{2}}. \n
\end{align}
The identities from $\, d(d(b_1)) =0, \, d(d(b_2)) =0$ determine the
derivative of $\, b_{11}$ by \beq d b_{11}=-\frac{8}{3} \,b_{{1
1}}\omega_{{0 0}} + ( 2\,b_{{1}}a_{{2}}+3\,b_{{1}}a_{{1}} )
(\omega^{{1}} - \omega^2).\n \eeq Moreover, $\, d(d(b_{11}))=0$ is
an identity.

At this step, we invoke a version of  Cartan-K\"ahler theorem.
\begin{proposition}\label{231R0}
The structure equation for the nondegenerate isothermally asymptotic
Legendrian surfaces with flat asymptotic web is in involution with
the general solution depending on one  arbitrary function  of 1
variable.
\end{proposition}
\Pf From the analysis above, the   exterior derivative identity
$\, d(d(c_1)) =0$  gives 1 compatibility equation  while the remaining
independent derivative coefficients at this step are
$\, \{ \, c_{11}; \, b_{11}  \, \}$. By inspection,  the resulting structure
equation is in involution with  the last nonzero Cartan character
s$_1=1$.
$\sq$

\subsection{Asymptotically ruled}\label{sec24}
In this sub-section, we consider the class of Legendrian surfaces
for which the asymptotic Legendrian $\, \PP^2$-field defined at
3-adapted frame is constant along the corresponding asymptotic
foliation.
\begin{definition}
Let $\, M \hook \PP^5$ be a nondegenerate Legendrian surface.
Let $\, \mathcal{F}$ be an asymptotic foliation(one of the three)
defined at 2-adapted frame.
$\, \mathcal{F}$ is \emph{ruled} if  the corresponding
Legendrian $\, \PP^2$-field defined at 3-adapted frame is constant
along $\,  \mathcal{F}$.
\end{definition}
\noi
With an  abuse of terminology, we call the leaf-wise constant
$\, \PP^2$-field \emph{rulings} of the asymptotic foliation.

We wish to give a differential analysis for  Legendrian surfaces
with three, two, or one ruled asymptotic foliations  in turn. Recall
$\, \omega^3=-(\omega^1+\omega^2)$. From the structure equation
\eqref{2frame},
\begin{align}\label{22ruled}
d e_0 &\equiv e_1 \, \omega^1, \mod  \, e_0; \omega^2,    \\
d e_1 &\equiv f_2 \, \omega^1, \mod  \, e_0, \, e_1; \omega^2, \n   \\
d f_2 &\equiv  e_2 \, (a_4 \,  \omega^1), \mod  \, e_0, \, e_1, \,
f_2; \omega^2, \n
\end{align}
\begin{align}
d e_0 &\equiv e_2 \, \omega^2, \mod  \, e_0; \omega^1, \n  \\
d e_2 &\equiv f_1 \, \omega^2, \mod  \, e_0, \, e_2; \omega^1, \n   \\
d f_1 &\equiv  e_1 \, (a_1 \,  \omega^2), \mod  \, e_0, \, e_2, \,
f_1; \omega^1, \n
\end{align}
\begin{align}
d e_0 &\equiv (e_1-e_2)  \, \omega^1, \mod  \, e_0; \omega^3, \n  \\
d (e_1-e_2)  &\equiv (f_1+f_2) \, \omega^1, \mod  \, e_0, \, (e_1-e_2) ; \omega^3, \n   \\
d (f_1+f_2)  &\equiv  e_1 \, (-(a_1 +2\, a_3) + (2\,a_2+a_4))\,
\omega^1, \mod  \, e_0, \, (e_1-e_2), \, (f_1+f_2); \omega^3. \n
\end{align}

\subsubsection{Tri-ruled}\label{sec241}
This is the class of  surfaces for which all of the  three
asymptotic foliations are ruled. The following differential analysis
shows that the local moduli space of tri-ruled Legendrian surfaces
consists of two points.

Assume that  each of the three Legendrian $\, \PP^2$-fields
$\, [e_0 \w e_1 \w f_2], \,  [e_0 \w e_2 \w f_1]$, and
$\,  [e_0 \w (e_1-e_2) \w (f_1+f_2)] $ is constant along the asymptotic
foliations defined by $\, \omega^2=0, \, \omega^1=0$, and
$\, \omega^3=0$ respectively.
By  \eqref{22ruled}, this implies
\beq
a_4=0, \; a_1=0, \; a_2 = a_3 =a. \n
\eeq
A tri-ruled Legendrian
surface has flat asymptotic 3-web, and it is also  isothermally
asymptotic. The cubic differential $\, \Psi$ vanishes when
$\, a = 0$, and  a non-flat tri-ruled  surface is  distinct from
$\, \Psi$-null surfaces discussed in Section \ref{sec22}.

A differential  analysis  shows that a tri-ruled surface necessarily
has $\, b_i=0, \, c_1=c_2=a^2$(we omit the details).
Maurer-Cartan form  $\, \phi$ is reduced  to
\beq\label{23triphi}
\phi   =
\begin{pmatrix}
\omega_{{0 0}}& a  \omega^{{1}}& a \omega^{{2}} &  a^2
(\omega^{{1}}+\omega^{{2}})
&\cdot&\cdot\\\noalign{\medskip}\omega^{{1}}&\frac{1}{3}
\,\omega_{{0 0}}
&\cdot&\cdot&\cdot&a(\omega^{{1}}+\omega^{{2}})\\
\noalign{\medskip}\omega^{{2}}&\cdot& \frac{1}{3} \,\omega_{{0
0}}&\cdot &a (\omega^{{1}}+\omega^{{2}})&\cdot
\\\noalign{\medskip}\cdot&\cdot&\cdot&-\omega_{{0 0}}&-\omega^{{1}}&-\omega^{{2}}
\\\noalign{\medskip}\cdot&\omega^{{2}}&\omega^{{1}}+\omega^{{2}}
&-a  \omega^{{1}}&-\frac{1}{3} \,\omega_{{0 0}}&\cdot \\
\noalign{\medskip}\cdot &\omega^{{1}}+\omega^{{2}}&\omega^{{1}} &-a
\omega^{{2}}&0&-\frac{1}{3} \,\omega_{{0 0}}
\end{pmatrix},
\eeq
with
\beq\label{24da}
d a = -\frac{4}{3} a \, \omega_{00}.
\eeq
\begin{proposition}
Let $\, M \hook \PP^5$ be a nondegenerate, tri-ruled Legendrian
surface. $\, M$ is necessarily isothermally asymptotic  with flat
asymptotic web. The Maurer-Cartan form of the 5-adapted frame of
$\, M$ is reduced to \eqref{23triphi}, and the single  structure
coefficient  $\, \{ \, a \, \}$ satisfies the   equation
\eqref{24da}. The local moduli space of  tri-ruled Legendrian
surfaces   consists of two points.
\end{proposition}
\Pf From \eqref{24da}, the moduli space is divided into two cases;
$\, a \equiv 0$, or $\, a \ne 0$. $\sq$

A differential geometric characterization of  tri-ruled  surfaces is
presented in   Section \ref{sec5}.

\subsubsection{Doubly-ruled}
Assume that  each of  the two  Legendrian $\, \PP^2$-fields
$\, [e_0 \w e_1 \w f_2]$ and  $\,  [e_0 \w e_2 \w f_1]$ is constant
along the asymptotic foliations  defined by $\, \omega^2=0$  and  $\omega^1=0$
respectively. By  \eqref{22ruled}, this implies
\beq\label{252ddouble} a_4=0, \; a_1=0. \eeq Note that such a doubly
ruled surface with flat asymptotic web is necessarily tri-ruled.

Successively differentiating \eqref{252ddouble}, one gets
\begin{align}
b_2&=0, \quad \quad   \quad \quad \, b_{3}=0, \n \\
a_{12}&=a_{41}=0, \quad    \,  b_{21}=b_{22}=b_{31}=0, \n \\
b_{{11}}&=-c_{{1}}+\frac{3}{5}\,{a_{{2}}}^{2}+\frac{2}{5}\,a_{{2}}a_{{3}}.
\n
\end{align}
Exterior derivative  $\, d(d(b_1)) =0$  with these relations  gives
\beq c_{21} ={\frac {44}{5}}\,b_{{1}}a_{{3}}-a_{{3}}a_{{2 1}}
+{\frac {16}{5}}\,b_{{1}}a_{{2}}-a_{{2}}a_{{2 1}}. \n \eeq The
identities from $\, d(d(a_2)) =0, \, d(d(a_3)) =0$ determine the
derivative of $\, a_{21}$ by
\begin{align}
d a_{21}&=-2\,a_{21} \omega_{00} + (-8\,c_{{2}}+{\frac
{48}{25}}\,{a_{{2}}}^{2}+{\frac {64}{25}}\,a_{{2}}a_{{3}} +{\frac
{88}{25}}\,{a_{{3}}}^{2})\omega^1 +(8\,c_{{1}}+12\,c_{{2}}+{\frac
{88}{25}}\,{a_{{2}}}^{2}-{\frac {56}{25}}\,a_{{2}}a_{{3}}
 -{\frac {132}{25}}\,{a_{{3}}}^{2})\omega^2. \n
\end{align}
Differentiating this equation again, one gets \beq c_{22}=c_{11}
-\frac{1}{2}\,a_{{2}}a_{{2 1}}+\frac{1}{2}\,a_{{3}}a_{{2 1}}
-8\,b_{{1}}a_{{3}}-2\,b_{{1}}a_{{2}}. \n \eeq

Exterior derivatives $\, d(d(c_1)) =0, \;d(d(c_2)) =0$  with these
relations determine the derivative of $\, c_{11}$(the exact
expression for $\, d c_{11}$ is long, and shall be omitted).
Differentiating this equation,   $\, d(d(c_{11})) =0$  finally gives
\begin{align}
c_{11}&= \frac{1}{4(a_{{2}}-a_{{3}})} (-15\,a_{{2
1}}c_{{1}}-7\,{a_{{3}}}^{2}a_{{2 1}}
+8\,a_{{2}}a_{{3}}b_{{1}}-32\,{a_{{2}}}^{2}b_{{1}}+80\,c_{{1}}b_{{1}}
+9\,{a_{{2}}}^{2}a_{{2 1}}-100\,c_{{2}}b_{{1}}\n \\
&\qquad  \qquad \qquad +44\,{a_{{3}}}^{2}b_{{1}}-2\,a_{{3}}a_{{2
1}}a_{{2}}+15\,c_{{2}}a_{{2 1}}). \n
\end{align}
Here we assumed  that  $\, a_2 \ne a_3$, or equivalently that the
Legendrian surface is not tri-ruled.

Differentiating this equation again, and comparing with the formula
for $\, d c_{11}$, one gets two polynomial compatibility equations
for   six   coefficients $\, \{ \, a_2, \, a_3, \, b_1, \, a_{21},
\, c_1, \, c_2 \, \}$. Successive derivatives of these   equations
generate  a sequence of    compatibility equations for a
nondegenerate Legendrian surface to  admit exactly two asymptotic
$\, \PP^2$-rulings.

Partly due to the complexity of the polynomial compatibility
equations, our analysis is incomplete. We suspect that if  there do
exist  nondegenerate, doubly-ruled(and not tri-ruled)  Legendrian
surfaces, the moduli space of such surfaces is at most discrete.

\subsubsection{Singly-ruled}
Assume  that   the  Legendrian $\, \PP^2$-field
$\,  [e_0 \w (e_1-e_2) \w (f_1+f_2)] $ is constant along the asymptotic curves
defined by $\, \omega^3=0$. From \eqref{22ruled}, this implies
\beq a_1+2\,a_3= 2\,a_2+a_4. \n \eeq Note the equivalence relations.
\begin{align}
\mbox{\emph{Singly ruled} and  \emph{flat asymptotic web} }&=
\mbox{\emph{Singly ruled} and  \emph{isothermally asymptotic}}, \n \\
&= \mbox{\emph{Isothermally asymptotic} and  \emph{flat asymptotic
web}}. \n
\end{align}

An analysis  shows that the structure equation for a nondegenerate,
singly-ruled Legendrian surface becomes  involutive  after  one
prolongation with the general solution depending  on  five
arbitrary functions of 1 variable. We omit the details of
differential analysis for  this case.


\section{Second order deformation}\label{sec3}
\begin{definition}\label{3deformdefi}
Let $\, \x: M \hook  \PP^5$  be a nondegenerate Legendrian surface.
Let $\, \x': M \hook  \PP^5$ be a   Legendrian deformation of
$\, \x$.
$\, \x'$ is a \emph{$k$-th order deformation} if there exists a
map $\, g:  M \to$ \textnormal{Sp}$_3\C$ such that for each
$\, p \in M$, the $\, k$-adapted frame bundles of $\, \x'$ and
$\, g(p)\circ \x$ are isomorphic at $\, p$. When the application map
$\, g$ is constant, the deformation $\, \x'$ is \emph{trivial}, and
$\, \x'$ is \emph{congruent} to $\, \x$  up to motion by
\textnormal{Sp}$_3\C$. Two deformations $\, \x_1'$ and $\, \x_2'$
are \emph{equivalent} if there exists an element $\, g_0 \in
\textnormal{Sp}_3\C$ such that $\, \x_2' = g_0 \circ \x_1'$.
A 'deformation'  would mean an 'equivalence class of deformations
modulo \textnormal{Sp}$_3\C$ action' for brevity.
\end{definition}
\noi It follows from the construction  of  adapted frames in
Section \ref{sec2} that a $\,  k$-th order deformation is a
$\, (k+1)$-th order deformation when for each $\, p \in M$, the
application map $\,  g(p)\in$ \textnormal{Sp}$_3\C$ not only
preserves the $\, k$-adapted frame at $\, p$, but also  the first
order derivatives of the $\, k$-adapted frame at $\, p$ (this is a
vague explanation, but the meaning is clear).

The definition of  $\, k$-th order  deformation indicates a way to
uniformize the various geometric conditions that naturally occur in
the theory of deformation and rigidity of submanifolds in a
homogeneous space. Take for an example the familiar case of surfaces
in  three dimensional Euclidean space with the usual 1-adapted
tangent frame of the group of Euclidean motions. One  surface  is a
first order deformation of  the other if  they have the same induced
metric, and it  is a  second  order deformation if  they also have
the same second fundamental form. By Bonnet's theorem, a  second
order deformation is  a congruence, \cite{Sp}.

Fubini, and Cartan studied the problem of third order  deformation
of  projective hypersurfaces in $\, \PP^{n+1}$, \cite{Ca2} and the
reference therein. For $\, n \geq 3$, a third order deformation of a
hypersurface with nondegenerate second fundamental form is
necessarily a congruence,  \cite{JM} for a modern  proof. For $\, n
= 2$, Cartan showed that a generic  surface does not admit a
nontrivial third order deformation, but that there exist two special
classes of surfaces which admit maximum three parameter family of
deformations.

The purpose of this section is to lay a foundation  for generalizing
Cartan's  work on projective deformation of surfaces in $\, \PP^3$
to deformation of  Legendrian surfaces in $\, \PP^5$. As a
preparation, we first consider the second order deformation. By
applying a modified moving frame method, we determine the
fundamental structure equation for  the second order deformation  of
a nondegenerate Legendrian surface. The analysis shows that the
resulting structure equation is in involution, and admits arbitrary
function worth solutions locally. This  implies that there is no
local obstruction to second order deformation of a Legendrian
surface.

The  structure equation established in this section will be applied
to the projective deformation of Legendrian surfaces with geometric
constraints in Section \ref{sec4}.

\subsection{Structure equation}
Let $\, \x: M \hook \PP^5$ be a nondegenerate  Legendrian surface.
Let F$\to M$ be the associated canonical bundle of  5-adapted frames
with the induced $\, \mathfrak{sp}_3 \C$-valued Maurer-Cartan form
$\, \phi$. The pair $\, (F, \, \phi)$ satisfies the properties
described in  Proposition \ref{2structure}.
Let $\, \x': M \hook \PP^5$ be a second order  deformation of  $\, \x$.
Let F$' \to M$ be
the associated canonical bundle with the induced Maurer-Cartan form
$\, \pi$. From the definition of   second order  deformation,
$\, \mbox{F}'$ can be considered as a graph over $\, \mbox{F}$ which
agrees with $\, \mbox{F}$ up to   2-adapted frame. By pulling back
$\, \pi$ on F, we  regard  $\, \pi$ as another
$\, \mathfrak{sp}_3 \C$-valued  Maurer-Cartan form on  F.

Set \beq\label{3setpi}
\pi = \phi + \del.
\eeq
The components of $\, \del$  are denoted by
\beq
\del =
\begin{pmatrix}
\delo & \dele \\
\delt  & - \delo^t
\end{pmatrix},\n
\eeq
where $\, \delt^t=\delt, \, \dele^t=\dele$.
Maurer-Cartan equations  for $\, \pi$ and $\, \phi$  imply the fundamental
structure equation for the deformation $\, \del$;
\beq\label{3fundast}
d(\del) + \del \w \phi + \phi \w \del + \del \w \del = 0.
\eeq
Differentiating the components of $\, \del$ from now
on would mean applying this structure equation.

We employ the method of moving frames to normalize the frame bundle
F$'$  based at F. In effect, one may adopt the  following analysis
as the constructive definition of  $(\,$F$', \pi \,)$. The
equivariant reduction process for  F$'$ in this section can be
considered as  the derivative of  the  one  applied for  F  in
Section \ref{sec2}. To avoid repetition, some of the details of
non-essential terms in the analysis below shall be omitted.

\vsp{0.1cm}
\texttt{1, and 2-adapted frame}. Let $\, (e', \, f')$ and $\, (e, \,
f)$ denote  the 5-adapted  Sp$_3 \C $-frames of  $\, \mbox{F}'$ and
$\,\mbox{F}$ respectively,   \eqref{2frame}. The condition of
second order deformation and the definition of 2-adapted frame imply
that there exist  frames such that
\begin{align}\label{3start}
e_0'  &=\, e_0,   \\
(  e_1', \, e_2' ) &\equiv (  e_1, \, e_2) \mod  e_0. \n
\end{align}
We take this identification as the initial circuit for the
algorithmic process of moving frame computation.

From the general theory of moving frames, \eqref{3start} shows that
one may adapt   F$'$  to normalize
\begin{align}\label{3stepfirsto}
\delo_{10} &=0, \, \delo_{20}=0,    \\
\delt_{00} &=0, \, \delt_{10}=0, \delt_{20}  =0.  \n
\end{align}
Differentiating these equations, one gets
\begin{align}
\begin{pmatrix}
\delo_{1 1}-\delo_{00} &\delo_{1 2} \\
\delo_{2 1}&\delo_{2 2}-\delo_{00}
\end{pmatrix} \w
\begin{pmatrix}
\omega^1 \\
\omega^2
\end{pmatrix} &=0,\n \\
\begin{pmatrix}
\delt_{1 1}&\delt_{1 2} \\
\delt_{2 1}&\delt_{2 2}
\end{pmatrix} \w
\begin{pmatrix}
\omega^1 \\
\omega^2
\end{pmatrix} &=0. \n
\end{align}
By Cartan's lemma, there exist coefficients $\, \delta s_{ijk}, \,
\delta t_{ijk};  \, i,j,k=1,2$,  such that
\begin{align}
\delo_{ij}- \delta_{ij} \delo_{00} &= \delta s_{ijk} \, \omega^k,
              \quad \; \; \mbox{where }   \delta s_{ijk} = \delta s_{ikj}, \n \\
\delt_{i j}  &=\delta t_{ijk} \, \omega^k,
              \quad \; \;  \mbox{where   $\,\delta t_{ijk}$ fully symmetric in indices}. \n
\end{align}

The coefficients $\,\{\,  \delta t_{ijk} \, \}$ depend  on the
second order jet of the immersion $\, \x'$. By the  assumption of
second order deformation, the cubic differential
\begin{align}
\Phi' &=(\theta_{ij} + \delt_{ij} ) \, (\omega^i+\delo_{i0}) ( \omega^j+\delo_{j0}),  \n \\
        &= \Phi + \delta t_{ijk}\, \omega^i \omega^j \omega^k, \quad \mbox{by} \; \eqref{3stepfirsto}, \n
\end{align}
must be a nonzero multiple of   $\, \Phi$. One may thus use the
group action that corresponds to $\, \delo_{00}$ to scale so that
\beq\label{3tijk}
\delta t_{ijk}=0.
\eeq

\texttt{3-adapted frame}. On the 2-adapted frame satisfying
\eqref{3tijk}, set $\, \delo_{00} = \delta s_{00k} \, \omega^k$.
There are $\, 6+2=8$ independent coefficients in
$\,\{ \, \delta s_{ijk}=\delta s_{ikj}, \, \delta s_{00k} \, \}$.
Differentiating
$\, \delt_{ij} = 0, \, i,j=1,2$, one gets  3 linear  relations among
them. By the group action that corresponds to
$\, \{ \delo_{01},\, \delo_{02}, \, \dele_{11}, \, \dele_{12}, \, \dele_{22} \}$,
one may translate the  remaining 5 coefficients so that
\begin{align}\label{3sijk}
\delo_{00}&=0,   \\
\delo_{ij}  &=0, \quad \mbox{for} \; \; i, j = 1, 2. \n
\end{align}
 At  this step,  the deformation $\, \del$ is reduced to
\beq\label{3delmatrix} \del =
\begin{pmatrix}
 \cdot &    \delo_{01}  &  \delo_{02}  &\dele_{{0 0}}&\dele_{{10}}&\dele_{{20}} \\
 \cdot &    \cdot           & \cdot           &  \dele_{{1 0}} &  \dele_{{11}}  &  \dele_{{12}}  \\
  \cdot&   \cdot            &   \cdot           &\dele_{{2 0}}&  \dele_{{12}}  &  \dele_{{22}} \\
 \cdot &   \cdot             &  \cdot            &  \cdot  &  \cdot   &   \cdot \\
  \cdot&   \cdot            &  \cdot               &  -  \delo_{01}   &   \cdot   &   \cdot \\
  \cdot&  \cdot            &  \cdot             &  -  \delo_{02}  &  \cdot   &   \cdot
\end{pmatrix}.
\eeq

Note that since all of the third  order terms
$\,\{ \, \delta s_{ijk}=\delta s_{ikj}, \, \delta s_{00k} \, \}$
are absorbed by
frame adaptations, \emph{a second order deformation of a
nondegenerate Legendrian surface is automatically a third order
deformation}, see remark below Definition \ref{3deformdefi}.

\texttt{4-adapted frame}. On the 3-adapted frame satisfying
\eqref{3sijk}, set
\begin{align}
\delo_{0i}&=\delta h_{0ik} \, \omega^k,  \quad \mbox{for} \; \; i=1, \, 2,   \n \\
\dele_{ij} = \dele_{ji}
                  &=\delta h_{ijk} \, \omega^k,
\quad \mbox{for} \; \; i, j = 1, 2. \n
\end{align}
There are   10 independent coefficients in  $\, \{ \, \delta
h_{0ij},  \, \delta h_{ijk} \, \}$. Differentiating   \eqref{3sijk},
one gets 5  linear  relations among them. By the group action that
corresponds to $\, \{  \dele_{10}, \,   \dele_{20} \}$, one may
translate $\, \delta h_{111}=0, \, \delta h_{222}=0$. The structure
coefficients can    be normalized  accordingly so that
\begin{align}\label{3eta111}
\delo_{01}&= (u_0+u_1)  \, \omega^1,   \quad
\delo_{02} =   (u_0+u_2) \, \omega^2,    \\
( \, \dele_{ij} \, ) &=
\begin{pmatrix}
u_2 \, \omega^2 &   u_0\, (\omega^1 +  \omega^2) \\
u_0\, (\omega^1 +  \omega^2) &  u_1 \, \omega^1
\end{pmatrix},  \n
\end{align}
for 3 coefficients $\, \{ u_0, \, u_1, \, u_2 \}$.

Note that the first equation of  \eqref{3eta111} implies that the
web curvature of the asymptotic 3-web is invariant under the second
order   deformation. Let us denote the covariant derivatives of
$\{ \, u_0, \, u_i \, \}$ by
\begin{align}
d u_0 &= -\frac{4}{3}\, u_0 \, \omega_{00} + u_{0k}\, \omega^k,  \n \\
d u_i  &= -\frac{4}{3}\,  u_i \, \omega_{00} + u_{ik}\,  \omega^k,
\quad \mbox{for} \; \, i=1, \, 2. \n
\end{align}

\texttt{5-adapted frame}. On the 4-adapted frame satisfying
\eqref{3eta111}, set \beq
 \dele_{i0}  = \delta h_{i0k} \, \omega^k, \quad \mbox{for} \; \; i = 1, 2. \n
\eeq By the group action that corresponds to $\, \{  \dele_{00} \,
\}$, one may translate $\, \delta h_{101}+\delta h_{202}=0$.
Introduce variables  $\,  \{  \,  v_0, \,v_1, \, v_2; \, w_1, \, w_2\, \}$,
and put
\begin{align}\label{3eta00}
 \dele_{10} &= v_0 \, \omega^1+ v_2 \, \omega^2,    \\
 \dele_{20} &= v_1 \, \omega^1 - v_0 \, \omega^2, \n   \\
 \dele_{00} &= w_1 \, \omega^1 + w_2 \, \omega^2. \n
\end{align}
At this step, no more frame adaptation is available. The reduction
process of  moving frame method  stops here. Let us denote the
covariant derivatives of   $\, \{ \,     v_i, \, w_i  \, \}$ by
\begin{align}
d v_i   &= - 2\,  v_i \, \omega_{00} + v_{ik}\,  \omega^k,  \n \\
d w_i  &= - \frac{8}{3} \,  w_i \, \omega_{00} + w_{ik}\,  \omega^k,
\quad \mbox{for} \; \, i=1, \, 2. \n
\end{align}

The normalization we chose  for the 4, and 5-adapted frame implies a
set of compatibility equations among the deformation coefficients
$\,  \{  \,  u_0, \, u_i; \, v_0, \,v_i; \, w_i\, \}$.
Differentiating $\, \{ \delo_{01}, \, \delo_{02} \}$ from
\eqref{3eta111}, one gets
\begin{align}\label{3u0}
u_{01} &=-u_{{2 1}}-2\,v_{{0}}+v_{{2}},     \\
u_{02} &= -u_{{1 2}}+2\,v_{{0}}+v_{{1}}. \n
\end{align}
Differentiating $\, \{ \dele_{11}, \,  \dele_{22},\,  \dele_{12}
\}$, one gets
\begin{align}\label{3v0}
v_0 &=\frac{1}{2}(-v_1+v_2),      \\
u_{12} &=-2\,v_{{1}}, \n \\
u_{21} &=-2\,v_{{2}}. \n
\end{align}
The identity from the exterior derivative $\, d(d(u_0))=0$  implies
\beq\label{3v12} v_{12}=2\,v_{{1 1}}+v_{{2 1}}-2\,v_{{2 2}}
+\frac{8}{5}\,(a_{{2}}    -  a_{{3}}) u_{{0}}. \eeq

Differentiating $\, \{ \dele_{10}, \, \dele_{20} \}$ from
\eqref{3eta00} with these relations, one gets
\begin{align}\label{3v11}
v_{11}&=3\,v_{{2 1}}-6\,w_{{1}}+8\,w_{{2}}-2\,{u_{{0}}}^{2} +{\frac
{16}{5}}\,a_{{2}} u_{{0}} -{\frac {36}{5}}\,a_{{3}} u_{{0}}
+2\,a_{{4}} u_{{2}}  +2\,a_{{1}}u_{{1}} +2\,u_{{1}}u_{{2}},   \\
v_{22}&=3\,v_{{2 1}}+6\,w_{{2}}-4\,w_{{1}}-2\,{u_{{0}}}^{2}+{\frac
{12}{5}}\,a_{{2}}u_{{0}} -{\frac {32}{5}}\,a_{{3}} u_{{0}}
+2\,a_{{4}} u_{{2}} +2\,a_{{1}}u_{{1}}+2\,u_{{1}}u_{{2}}. \n
\end{align}
Differentiating $\, \{ \dele_{00}  \}$ from \eqref{3eta00}, one
finally gets
\begin{align}\label{3w12}
w_{12}&=w_{{2 1}}+ ( 2\,v_{{2}}-2\,v_{{1}} +2\,b_{{2}} -2\,b_{{3}}
) u_{{0}} +  ( 2\,v_{{2}}+2\,b_{{2}}  ) u_{{1}}
+ ( -2\,v_{{1}}-2\,b_{{3}}  ) u_{{2}}   \\
&+  ( -2\,a_{{1}}-2\,a_{{3}}  ) v_{{1}}+
        ( 2\,a_{{4}}+2\,a_{{2}}  ) v_{{2}}.  \n
\end{align}
The fundamental structure equation \eqref{3fundast} for $\, \del$ is
now an identity.

\begin{proposition}\label{3deform2}
Let $\, \x: M \hook \PP^5$ be a nondegenerate   Legendrian surface.
Let $\,  \x' : M \hook \PP^5$ be  a second order   Legendrian
deformation of $\, \x$. Let $\, \pi = \phi+ \del$ be the induced
Maurer-Cartan form of   $\, \x'$, \eqref{3setpi}, where $\, \phi$ is
the  induced Maurer-Cartan form of   $\, \x$.
There exists a 5-adapted frame  for $\, \x'$ such that the coefficients of
$\, \del$   satisfy  the structure equations \eqref{3delmatrix} through
\eqref{3w12}. These   equations   furthermore imply that;

\textnormal{a)} The  structure equation  for second order
deformation becomes  involutive after one prolongation with the
general solution depending on five arbitrary functions of 1
variable.

\textnormal{b)} The second order   deformation $\, \x'$ is
necessarily a third order deformation. If    $\, \x'$ is  a fourth
order deformation of  \x, $\, \x'$ is congruent to $\, \x$.
\end{proposition}
\Pf a) We show that the structure equation  for deformation becomes
involutive after a partial prolongation. The identities from
exterior derivatives $\, d(d(v_1))=0, \, d(d(v_2))=0$ determine the
derivative of $\, v_{21}$ by \beq\label{3v21} d v_{21}\equiv
(-w_{11}-3\,w_{22}+3\,w_{21})\omega^1
                       + ( w_{11}- \,w_{22}-\,w_{21})\omega^2,
                       \mod \; \omega_{00}; \,u_{11}, \, u_{22}; \, u_0, \, u_i, \, v_i, \, w_i.
\eeq Note the relation
\begin{align}\label{3rest}
\begin{pmatrix}
d u_1 \\
d u_2 \\
d w_1 \\
d w_2
\end{pmatrix} \equiv
\begin{pmatrix}
u_{11} & \cdot \\
\cdot & u_{22} \\
w_{11} & w_{21} \\
w_{21} & w_{22}
\end{pmatrix}
\begin{pmatrix}
\omega^1 \\
\omega^2
\end{pmatrix}
\mod \;\omega_{00}; \,u_0, \, u_i, \, v_i, \, w_i.
\end{align}
By inspection, the structure equations \eqref{3v21} and
\eqref{3rest} are in involution with the last nonzero Cartan
character  s$_1=5$.

b) For the first part, see the remark at the end of
\texttt{3-adapted frame}. For the second part, the condition for the
fourth order deformation implies $\, u_0=u_1=u_2=0$. The
compatibility equations  \eqref{3v0} and \eqref{3v11} then force
the remaining deformation coefficients to vanish so that
$\, \del=0$. The rest follows from the uniqueness theorem of ODE,
\cite{Gr}. $\sq$

We shall examine the second order  Legendrian deformation with the
additional condition that it preserves the fourth order differential
$\, \Psi$, or that it preserves both $\, \Psi$ and the fifth order
differential $\, \chi$. The primary object of  our  analysis will be
to give  characterization of  such  surfaces that support maximum
parameter family of  nontrivial deformations.

\section{$\Psi$-deformation}\label{sec4}
In this section, we apply the fundamental structure equation for
second order deformation to  the geometric situation where the
deformation is required to preserve a part of  fourth order
invariants of a Legendrian surface.
\begin{definition}\label{4deformdefi}
Let $\, \x: M \hook  \PP^5$  be a nondegenerate Legendrian surface.
Let $\, \x': M \hook  \PP^5$ be a second order  deformation of
$\, \x$.
$\, \x'$ is a \emph{$\, \Psi$-deformation} if the application
map $\, g:  M \to$ \textnormal{Sp}$_3\C$ for the second order
deformation $\, \x'$ is such that for each $\, p \in M$, the fourth
order cubic differential $\, \Psi'$ of  $\, \x'$ and $\, \Psi$ of
$\, g(p)\circ \x$ are isomorphic at $\, p$.
\end{definition}
As noted in  Proposition \ref{3deform2}, there is no local
obstruction  for the  second order deformation of a Legendrian
surface, whereas if one requires the second order deformation to
preserve all of the fourth order  invariants, the deformation is
necessarily a congruence. The idea is to  impose a condition that
balances between these two extremes.

Let us give a summary of   results in this section. The condition
for a second order deformation to be a $\, \Psi$-deformation is
expressed as  a pair of linear equations on the deformation
coefficients, \eqref{4delcon}. A more or less  basic over-determined
PDE analysis of these equations shows that the resulting  structure
equation for  $\, \Psi$-deformation closes up admitting at most
three parameter family of solutions, \eqref{4delstructure}.
The class of  isothermally asymptotic surfaces with flat asymptotic web
discussed in Section \ref{sec231} are examples of  such  surfaces
admitting maximum  parameter  family of  $\, \Psi$-deformations,
which we call \emph{$D_0$-surfaces}, \eqref{R0defi}.
Analysis of the structure
equation shows that there exist another class of surfaces with
finite local moduli that admit maximum parameter family of  $\, \Psi$-deformations,
which we call \emph{$D$-surfaces}, \eqref{Rdefi}. $D_0$-surfaces and
$D$-surfaces account for the set of maximally $\, \Psi$-deformable
Legendrian surfaces,  Theorem \ref{4main}. Further analysis shows
that there exist subsets called  $S_0$-surfaces and $S$-surfaces
which admit $\, \Psi$-deformations that also preserve the fifth
order  differential $\, \chi$.

We continue the  analysis of  Section \ref{sec3}.

\subsection{Structure equation}\label{sec41}
Let $\, \x: M \hook \PP^5$ be a nondegenerate  Legendrian surface.
Let $\,  \x' : M \hook \PP^5$ be  a  $\, \Psi$-deformation of $\, \x$.
Let $\, \pi = \phi+ \del$ be the induced Maurer-Cartan form of
$\, \x'$, where $\, \phi$ is  the induced Maurer-Cartan form of $\,\x$.
From \eqref{3eta111}, the deformation of the invariant
differentials $\, \Psi$ and $\, \chi$ are given by
\begin{align}\label{4delPsi}
\delta \Psi &= \delta (\eta_{ij} \omega^i  \omega^j ),   \\
                 &= (u_2+2 \, u_0)\, (\omega^1)^2   \omega^2 + (u_1+2 \, u_0)\,\omega^1  ( \omega^2)^2, \n \\
\delta \chi &= \delta (\eta_{10} \omega^1+\eta_{20} \omega^2), \n \\
                 &= v_0\, (\omega^1)^2+ (v_1+v_2) \omega^1  \omega^2 - v_0\, (\omega^2)^2. \n
\end{align}
The condition for the deformation to preserve $\, \Psi$ is expressed
by the pair of linear equations
\beq\label{4delcon}
u_1+2 \, u_0=0,
\; \; u_2+2 \, u_0=0.
\eeq
We wish to  give an  analysis of the
compatibility equations for  the deformation $\, \del$ derived from
\eqref{4delcon}.

Differentiating \eqref{4delcon}, one gets
\begin{align}\label{4delv}
v_2&=-v_1,   \\
u_{11}&=  2\,v_1, \n \\
u_{22}&= -2\,v_1. \n
\end{align}
Since $\, v_0 = \frac{1}{2}(-v_1+v_2)$, we observe that
\emph{a $\, \Psi$-deformation leaves $\, \chi$ invariant when $\, v_1=v_2=0$}.

Differentiating $\, v_1+v_2=0$,  one gets
\begin{align}\label{4delw}
w_2     &= w_{{1}}+  \frac{6}{5} ( -\,a_{{2}}+\,a_{{3}}) u_{{0}},   \\
v_{21}&=-\frac{1}{2}\,w_{{1}}-\frac{3}{2}\,{u_{{0}}}^{2} + (
a_{{4}}+a_{{1}}+\frac{8}{5}\,a_{{2}}-\frac{3}{5}\,a_{{3}}) u_{{0}}.
\n
\end{align}
Differentiating the first equation of \eqref{4delw} for  $\, w_2$,
one gets
\begin{align}
w_{21}&= w_{{1 1}}+\frac{6}{5} (  \,a_{{2}}- \,a_{{3}} ) v_{{1}}
+  ( 9\,b_{{3}}+ 12\,b_{{1}}-3\,a_{{2 1}}  ) u_{{0}}, \n \\
w_{22}&=w_{11}+4\,v_{{1}}u_0+(
-2\,a_{{4}}-2\,a_{{1}}-2\,a_{{3}}-2\,a_{{2}}) v_{{1}} + (
-12\,b_{{1}}-4\,b_{{3}}+4\,b_{{2}}) u_{{0}}. \n
\end{align}
The identity from the exterior derivative   $\, d(d(v_1))=0$ with
these relations finally gives \beq\label{4w11} w_{11}
=-2\,u_{{0}}v_1+( -{\frac {12}{5}}\,a_{{3}}+{\frac {22}{5}}\,a_{{2}}
+a_{{4}}+a_{{1}} ) v_{{1}}+  ( 24\,b_{{1}}-4\,b_{{2}}+15\,b_{{3}}
+a_{{4 1}}+ a_{{1 2}}-3\,a_{{2 1}} ) u_{{0}}. \eeq

At this step, the remaining independent deformation coefficients are
$\, \{ \, u_0, \, v_1, \, w_1 \}$. Moreover, they satisfy a closed
structure equation, i.e., their derivatives are expressed as
functions of  themselves and do  not involve any new variables. Let
us record the structure equations for
$\, \{ \, u_0, \, v_1, \, w_1 \}$.
\begin{align}\label{4delstructure}
d u_0 +\frac{4}{3}\,u_{{0}}\omega_{{00}}&
= v_{{1}}(- \omega^{{1}}+ \omega^{{2}}),   \\
d v_1 +2\,v_{{1}}\omega_{{00}}& =(\frac{3}{2} \,{u_{{0}}}^{2} + (
\frac{3}{5}\,a_{{3}}-a_{{4}}-a_{{1}}-\frac{8}{5}\,a_{{2}}) u_{{0}}
+\frac{1}{2}\,w_{{1}}) \omega^1\n \\
&\quad +(-\frac{3}{2} \,{u_{{0}}}^{2}+ ( a_{{3}}+a_{{4}}+a_{{1}})
u_{{0}}
-\frac{1}{2}\,w_{{1}})\omega^2, \n \\
d w_1 +\frac{8}{3}\,w_{{1}}\omega_{{00}}& = w_{11} \omega^1+ w_{12}
\omega^2, \n
\end{align}
where $\, w_{11}$ is in \eqref{4w11} and \beq w_{12} =  2\,u_0
v_{{1}} + ( -{\frac {28}{5}}\,a_{{3}} +{\frac
{18}{5}}\,a_{{2}}-a_{{4}}-a_{{1}} ) v_{{1}} +(
36\,b_{{1}}-6\,b_{{2}}+26\,b_{{3}}+a_{{4 1}}+a_{{1 2}}-6\,a_{{2 1}})
u_{{0}}. \n \eeq

Exterior derivative   $\, d(d(w_1))=0$ gives a universal
integrability condition  for $\, \Psi$-deformation;
\beq\label{4universal}
 (  -60\,b_{{1}}+8\,b_{{2}}  -42\,b_{{3}}-a_{{4 1}}-a_{{1 2}} +10\,a_{{2 1}} ) v_{{1}}
 + {\frac {12}{5}} ( \,a_{{2}}- \,a_{{3}}  ) w_{{1}}
 \equiv 0 \mod \; u_0.
\eeq The full expression for the right hand side of
\eqref{4universal} is given by
\begin{align}
\textnormal{RHS of  \eqref{4universal}} &=
- \frac {44}{15}\,( a_{{2}}- a_{{3}}) u_{{0}}^2
+(-5\,c_{{1}}+8\,c_{{2}}-17\,b_{{11}}+10\,b_{{2 1}}-\frac{4}{3}\,b_{{2 2}}
-{\frac {26}{3}}\,b_{{31}} \n \\
&\quad \; \; +{\frac {32}{15}}\,a_{{4}}a_{{2}}-{\frac
{32}{15}}\,a_{{4}}a_{{3}} -{\frac {112}{25}}\,{a_{{3}}}^{2}+{\frac
{123}{25}}\,{a_{{2}}}^{2}
-{\frac {86}{25}}\,a_{{2}}a_{{3}}-{\frac {32}{15}}\,a_{{1}}a_{{3}}\n \\
&\quad \; \; +{\frac {32}{15}}\,a_{{1}}a_{{2}}+3\,a_{{4}}a_{{1}}
-\frac{1}{3}\,a_{{1 2 1}}+2\,a_{{2 1 1}} -\frac{1}{3}\,a_{{4 1
1}}-a_{{2 1 2}} +\frac{1}{3}\,a_{{1 2 2}}+\frac{1}{3}\,a_{{4 1 2}})
u_0,  \n
\end{align}
where $\, a_{ijk}$ denote the covariant derivative  of $\, a_{ij}$
as before.

\begin{proposition}\label{41pro}
Let $\, M \hook \PP^5$ be a nondegenerate  Legendrian surface.
A $\,\Psi$-deformation of  $\, M$ is determined by three parameters
$\, \{ \, u_0, \, v_1, \, w_1 \, \}$ by  \eqref{4delcon}, \eqref{4delv},
and \eqref{4delw}. These three deformation parameters satisfy a
closed structure equation \eqref{4delstructure}. The structure
equation, and the universal integrability condition
\eqref{4universal} imply that;

\textnormal{a)} A nondegenerate Legendrian surface admits at most
three parameter family of $\, \Psi$-deformations.

\textnormal{b)} If the asymptotic web of   $\, M$ is not flat,
$\, M$ admits at most   two  parameter family of $\, \Psi$-deformations.

\textnormal{c)} Assume the asymptotic web of $\, M$ is flat. If  the
structure coefficients of $\, M$ do  not satisfy the differential
relation $\,( a_{{4 1}}+a_{{1 2}}+4\,b_{{2}}-4\,b_{{1}})=0$, $\, M$
admits at most one  parameter family of $\, \Psi$-deformations.

\textnormal{d)} Assume  the asymptotic web of   $\, M$ is flat and
the structure coefficients satisfy the relation
$\,( a_{{4 1}}+a_{{1 2}}+4\,b_{{2}}-4\,b_{{1}})=0$. If the  structure
coefficients of   $\, M$ do not  satisfy the  additional relation
$\, c_1=c_2$,  $\, M$ does not admit nontrivial
$\, \Psi$-deformations.

\textnormal{e)} A nondegenerate Legendrian surface $\, M$  admits
maximum three parameter family of $\, \Psi$-deformations if, and
only if $\, M$ has flat asymptotic web, and the structure
coefficients of $\, M$ satisfy the following   relations.
\begin{align}
a_2&=a_3, \n \\
b_1-b_2&= \frac{1}{4} (a_{{4 1}}+a_{{1 2}}), \n \\
c_1 &=c_2. \n
\end{align}
\end{proposition}
\Pf a) It follows from the uniqueness theorem of ODE, \cite{Gr}.

b) If the web curvature  \eqref{21webK} of the asymptotic web does
not vanish identically, one can solve \eqref{4universal} for
$\, w_1$ on a dense open subset of $\, M$.

c) The  asymptotic web is flat when  $\, a_2-a_3=0$. By the
structure equations from Section \ref{sec21}, \eqref{4universal} is
reduced to
\beq\label{4universal2}
v_{{1}} ( a_{{4 1}}+a_{{1 2}}+4\,b_{{2}}-4\,b_{{1}} ) \equiv 0  \mod \; \; u_0.
\n
\eeq
Under the assumption of c), one can solve for $\, v_1$ as a function of
$\, u_0$. Differentiating this,  \eqref{4delstructure} implies that
$\, w_1$ is also determined as a function of $\, u_0$.

d) and e) When $\,  a_2-a_3=0$  and
$\,  a_{{4 1}}+a_{{1 2}}+4\,b_{{2}}-4\,b_{{1}} =0$,
\eqref{4universal} is reduced to
\beq\label{4universal3}
u_0 (c_1-c_2)=0.\n
\eeq
If $\, c_1-c_2 = 0$,
the structure equation  \eqref{4delstructure} is compatible and
admits solutions with maximum three dimensional moduli.
If $\, c_1-c_2$ does not vanish identically,  $\, u_0=0$.
The structure
equation \eqref{4delstructure} then  implies $\, v_1=w_1=0$.
$\sq$

\begin{example}\label{threeexample}
The  analysis of  Section \ref{sec2} shows  that the following
classes of  Legendrian surfaces admit maximum three parameter family
of  nontrivial $\, \Psi$-deformations.

\emph{a)  $\Psi$-null surfaces, Section \ref{sec22}}

\emph{b)  Isothermally asymptotic surfaces with flat asymptotic web,
Section \ref{sec231}}

\emph{c)  Tri-ruled surfaces, Section \ref{sec241}}

\noi Note that  \emph{a)} and  \emph{c)} are subsets of   \emph{b)}.
\end{example}


It is evident that a generic  nondegenerate Legendrian surface does
not admit any nontrivial $\, \Psi$-deformations. In consideration of
the main theme of the paper, to understand Legendrian surfaces  with
special characteristics, we do not pursue to formulate the explicit
criteria for   $\, \Psi$-rigidity.

\subsection{Surfaces with maximum $\, \infty^{3}$
$\, \Psi$-deformations}\label{sec42}
The structure of the moduli space of solutions to the deformation
equation \eqref{4delstructure} depends on the geometry of the base
Legendrian surface. Among the variety of cases, we consider in this
subsection the class of Legendrian surfaces that admit maximum three
parameter family of $\, \Psi$-deformations. The rationale for this
choice comes from the fact that Kummer's quartic surface constitutes
an example of Cartan's maximally third order deformable surfaces in
$\, \PP^3$, \cite{Fe}.

Let $\, M \hook \PP^5$ be a nondegenerate  Legendrian surface with
maximum three parameter family of $\, \Psi$-deformations. From e) of
Proposition \ref{41pro}, such surfaces are characterized by the
following three relations on the  structure coefficients;
\begin{align}\label{42defining}
a_2-a_3&=0,   \\
b_1-b_2&= \frac{1}{4} (a_{{4 1}}+a_{{1 2}}), \n \\
c_1-c_2&=0. \n
\end{align}
We wish to  give an  analysis of the compatibility conditions
derived from these relations, and  determine the  structure equation
for the maximally $\, \Psi$-deformable surfaces.

Since $\, M$ has flat asymptotic web, let us assume the results of
Section \ref{sec21} and continue the analysis  from that point on.
Differentiating the third equation of \eqref{42defining}, one gets
\begin{align}\label{4201}
c_{21}&=c_{11},   \\
c_{22}&=c_{11}+ ( 4\,a_{{1}}+4\,a_{{2}}) b_{{1}}+ (
-2\,a_{{1}}+2\,a_{{4}} ) b_{{2}}. \n
\end{align}
The remaining undetermined derivative coefficients at this step are
$\, \{ \, a_{41},  c_{11}; \, b_{11} \}$. The identities from exterior
derivatives $\, d(d(a_1))=0, \, d(d(a_4))=0$ determine the
derivative of  $\, a_{41}$ by \beq\label{4202}
da_{41}+2\,a_{41}\omega_{00}= 2\, ( b_{{1 1}}
-{a_{{2}}}^{2}+c_{{1}}+a_{{4}}a_{{1}} ) ( \omega^{{1}}+ \omega^{{2}}
)  +4\,b_{11}\omega^2. \eeq Moreover, $\, d(d(a_{41}))=0$ is an
identity.

Exterior derivative $\,  d(d(b_{11}))=0$, \eqref{21b11}, with these
relations gives the universal integrability condition to admit
maximum $\, \infty^{3}$  family of  $\, \Psi$-deformations.
\beq\label{42uni} ( a_{{1}}-a_{{4}} ) c_{{1}}+ 2 ( a_{{1}}-
\,a_{{4}} ) b_{{1 1}} +3( -\,b_{{1}} + \,b_{{2}} ) (a_{{4 1}}+
2\,b_2) + ( a_{{1}}-a_{{4}} )  ( a_{{4}}a_{{1}}-{a_{{2}}}^{2})=0.
\eeq \noi At this juncture, the analysis divides into two cases.

\texttt{Case  $\, a_1-a_4=0$}. This is the case  of isothermally
asymptotic surfaces with flat asymptotic web. As noted in Example
\ref{threeexample}, this class of surfaces satisfy the defining
relations \eqref{42defining} and admit maximum three parameter
family of $\, \Psi$-deformations.

\begin{definition}\label{R0defi}
A \emph{Legendrian $D_0$-surface} is an immersed, nondegenerate
Legendrian surface in $\, \PP^5$ which is isothermally asymptotic
with flat asymptotic 3-web.
\end{definition}

Let us record the full structure equation  for $D_0$-surfaces.
\begin{align}\label{R0st}
da_1+\frac{4}{3} \,a_{{1}}\omega_{{0 0}}&= -2\,b_{{2}} \omega^{{1}}
+  ( 4\,b_{{1}}-2\,b_{{2}}   ) \omega^{{2}},   \\
da_2+\frac{4}{3}\,a_{{2}}\omega_{{0 0}}&=
 ( -2\,b_{{1}}+3\,b_{{2}} ) \omega^{{1}}
+ ( -4\,b_{{1}}+3\,b_{{2}}  ) \omega^{{2}},  \n \\
db_1+2\,b_{{1}}\omega_{{0 0}}&=
 b_{{1 1}}( \omega^{{1}}- \omega^{{2}}), \n \\
db_2+2\,b_{{2}}\omega_{{0 0}}&=
  ( -b_{{1 1}}+{a_{{2}}}^{2}-c_{{1}}-a_{{1}}^2)( \omega^{{1}}+\omega^2)
-2 \,b_{{1 1}} \, \omega^{{2}}, \n \\
dc_{1}+\frac{8}{3}\,c_{{1}}\omega_{00}&= c_{{1 1}}
(\omega^{{1}}+\omega^2)
+ 4\,b_{{1}}  (  a_{{1}}+ a_{{2}} )  \omega^{{2}},  \n \\
db_{11}+\frac{8}{3} \,b_{{11}}\omega_{{0 0}}&= b_{{1}}  (
2\,a_{{2}}+3\,a_{{1}} )   (  \omega^{{1}}- \omega^{{2}} ). \n
\end{align}
\noi The induced Maurer-Cartan form $\, \phi$ takes the following
form.
\beq \phi =
\begin{pmatrix}
\omega_{{0 0}}& (a_1+a_2) \omega^1&(a_1+a_2) \omega^2
&c_{{1}}(\omega^1+\omega^2) &b_{{1}}\omega^1+b_{{2}}\omega^2
&(-2\,b_{{1}} +b_{{2}}) \omega^1-b_{{1}}\omega^2
\\\noalign{\medskip}\omega^1&\frac{1}{3} \,\omega_{{0 0}}&\cdot
&b_{{1}}\omega^1+b_{{2}}\omega^2&a_{{1}}\omega^2
&a_{{2}}(\omega^1+\omega^2)
\\\noalign{\medskip}\omega^2&\cdot&\frac{1}{3}\,\omega_{{0 0}}
&(-2\,b_{{1}}+b_{{2}}) \omega^1-b_{{1}}\omega^2
&a_{{2}}(\omega^1+\omega^2)&a_{{1}}\omega^1
\\\noalign{\medskip}\cdot&\cdot&\cdot&-\omega_{{0 0}}&-\omega^1&-\omega^2
\\\noalign{\medskip}\cdot&\omega^2&\omega^1+\omega^2
&-(a_1+a_2) \omega^1&-\frac{1}{3} \,\omega_{{0 0}}&\cdot
\\\noalign{\medskip}\cdot&\omega^1+\omega^2&\omega^1
&-(a_{{1}}+a_2) \omega^2&\cdot&-\frac{1}{3} \,\omega_{{0 0}}
\end{pmatrix}
\eeq

\texttt{Case  $\, a_1-a_4 \ne 0$}. The structure equation closes up
in this case. First,  solve \eqref{42uni} for $\, c_1$.
Differentiating this, one can solve for $\, c_{11}$. At this step,
the structure equation for this class of surfaces closes up with 7
independent  structure coefficients $\, \{ \, a_1, \, a_2, \, a_4,
\, b_1, \, b_2; \, a_{41}, \, b_{11} \, \}$. Moreover,  an analysis
shows that the resulting structure equation is  compatible, i,e,,
$\, d^2=0$ is an identity and does not impose any new compatibility
conditions.

\begin{definition}\label{Rdefi}
A \emph{Legendrian $D$-surface} is an immersed, nondegenerate
Legendrian surface in $\, \PP^5$ which satisfies the following
conditions.

a) it is not isothermally asymptotic,

b) it has flat asymptotic 3-web, and satisfies the differential
relations \eqref{42defining}, \eqref{4201}, \eqref{4202}, and
\eqref{42uni}.
\end{definition}

Let us record the full structure equation  for $D$-surfaces.
\begin{align}\label{Rst}
da_1+\frac{4}{3} \,a_{{1}}\omega_{{0 0}}&= -2\,b_{{2}} \omega^{{1}}
+  ( -a_{{4 1}}-4\,b_{{2}}+4\,b_{{1}}  ) \omega^{{2}},   \\
da_2+\frac{4}{3}\,a_{{2}}\omega_{{0 0}}&=
 ( -2\,b_{{1}}+3\,b_{{2}}  ) \omega^{{1}}
+ ( -4\,b_{{1}}+3\,b_{{2}}  ) \omega^{{2}},  \n \\
da_4+\frac{4}{3}\,a_{{4}}\omega_{{0 0}}&=
 a_{{4 1}} \omega^{{1}}
+ ( 4\,b_{{1}}-2\,b_{{2}} ) \omega^{{2}},  \n \\
db_1+2\,b_{{1}}\omega_{{0 0}}&=
 b_{{1 1}}( \omega^{{1}}- \omega^{{2}}), \n \\
db_2+2\,b_{{2}}\omega_{{0 0}}&= - (  b_{{1
1}}-{a_{{2}}}^{2}+c_{{1}}+a_{{4}}a_{{1}})( \omega^{{1}}+\omega^2)
-2 \,b_{{1 1}} \, \omega^{{2}}, \n \\
da_{41}+2\,a_{41}\omega_{00} &=
 2\, ( b_{{1 1}} -{a_{{2}}}^{2}+c_{{1}}+a_{{4}}a_{{1}} )
 ( \omega^{{1}}+ \omega^{{2}}  )  +4\,b_{11}\omega^2,  \n \\
db_{11}+\frac{8}{3} \,b_{{11}}\omega_{{0 0}}&= -\frac{1}{4} \,  (
-8\,b_{{1}}a_{{1}}+4\,b_{{2}}a_{{1}} +a_{{4
1}}a_{{1}}-8\,b_{{1}}a_{{2}}-4\,b_{{1}}a_{{4}}+a_{{4 1}}a_{{4}} )
 ( \omega^{{1}}- \omega^{{2}} ). \n
\end{align}
The induced Maurer-Cartan form $\, \phi$ takes the following form.
\beq \phi =
\begin{pmatrix}
\omega_{{0 0}}& (a_4+a_2) \omega^1&(a_1+a_2) \omega^2
&c_{{1}}(\omega^1+\omega^2) &b_{{1}}\omega^1+b_{{2}}\omega^2
&(-2\,b_{{1}} +b_{{2}}) \omega^1-b_{{1}}\omega^2
\\\noalign{\medskip}\omega^1&\frac{1}{3} \,\omega_{{0 0}}&\cdot
&b_{{1}}\omega^1+b_{{2}}\omega^2&a_{{1}}\omega^2
&a_{{2}}(\omega^1+\omega^2)
\\\noalign{\medskip}\omega^2&\cdot&\frac{1}{3}\,\omega_{{0 0}}
&(-2\,b_{{1}}+b_{{2}}) \omega^1-b_{{1}}\omega^2
&a_{{2}}(\omega^1+\omega^2)&a_{{4}}\omega^1
\\\noalign{\medskip}\cdot&\cdot&\cdot&-\omega_{{0 0}}&-\omega^1&-\omega^2
\\\noalign{\medskip}\cdot&\omega^2&\omega^1+\omega^2
&-(a_4+a_2) \omega^1&-\frac{1}{3} \,\omega_{{0 0}}&\cdot
\\\noalign{\medskip}\cdot&\omega^1+\omega^2&\omega^1
&-(a_{{1}}+a_2) \omega^2&\cdot&-\frac{1}{3} \,\omega_{{0 0}}
\end{pmatrix},
\eeq where $\, c_1$ is given by \eqref{42uni}.

\begin{theorem}\label{4main}
The set of nondegenerate Legendrian surfaces in $\, \PP^5$ which
admit  maximum three parameter family of $\, \Psi$-deformations fall
into two categories; Legendrian $\, D_0$-surfaces, or Legendrian
$\,D$-surfaces. A general Legendrian $\, D_0$-surface depends on one
arbitrary function of 1 variable, whereas a  general Legendrian
$\,D$-surface depends on four constants.
\end{theorem}
\Pf The generality of solutions for the structure equation for
$\,D_0$-surfaces is treated in Proposition \ref{231R0}. For the
generality of  $\, D$-surfaces, consider  the invariant map
$\,I=(a_1, \, a_2, \, a_4, \, b_1, \, b_2; \, a_{41}, \, b_{11} ):$
F$\to \C^7$, where F   is the canonical bundle of 5-adapted frames.
Since $\, I$ generically has rank 3, the local moduli space of
$\,D$-surfaces has general dimension dim $(\C^7) -$ rank$(I) =4$.
$\sq$

The analogy of  Theorem \ref{4main} with Cartan's classification of
maximally third order deformable surfaces in $\, \PP^3$ is obvious,
\cite{Ca2}. Cartan's classification is also divided into two cases;
one case with infinite dimensional local moduli,
and the other case  with finite dimensional local moduli.
This analogy in a way conversely
justifies our choice of $\, \Psi$-deformations.

\subsection{$(\Psi, \, \chi)$-deformations}\label{sec43}
In this subsection, we examine which of the maximally
$\,\Psi$-deformable surfaces admit deformations that  leave  invariant
both $\, \Psi$  and the fifth order quadratic differential
$\,\chi$, \eqref{2Psi}.
\begin{definition}\label{5deformdefi}
Let $\, \x: M \hook  \PP^5$  be a nondegenerate Legendrian surface.
Let $\, \x': M \hook  \PP^5$ be a $\, \Psi$-deformation of  $\, \x$.
$\, \x'$ is a \emph{$(\Psi, \, \chi)$-deformation} if the
application  map $\, g:  M \to$ \textnormal{Sp}$_3\C$ for
the $\, \Psi$-deformation $\, \x'$ is such that for each $\, p \in M$,
the fifth order quadratic differentials $\, \chi'$ of  $\, \x'$ and
$\, \chi$ of  $\, g(p)\circ \x$ are isomorphic at $\, p$.
\end{definition}

Let us give a summary of   results in this subsection. The condition
for a  $\, \Psi$-deformation to be a $(\Psi, \, \chi)$-deformation
is expressed  by a single  linear equation on the deformation
coefficients, \eqref{43delcon}. An over-determined PDE analysis of
this equation shows that the resulting  structure equation for
$(\Psi, \, \chi)$-deformation closes up admitting at most one
parameter family of solutions, \eqref{43delstructure}. The structure
equation for the subset of maximally $\, \Psi$-deformable surfaces
which admit  one parameter family of
$(\Psi, \, \chi)$-deformations is then determined, Theorem \ref{43thm}.

We continue the  analysis of  Section \ref{sec42}, specifically from
\eqref{42uni}.

Let $\, \x: M \hook \PP^5$ be a nondegenerate, maximally
$\, \Psi$-deformable Legendrian surface. Let $\,  \x' : M \hook \PP^5$
be  a  $\, \Psi$-deformation of $\, \x$. Let $\, \pi = \phi+ \del$
be the induced Maurer-Cartan form of   $\, \x'$, where $\, \phi$ is
the induced Maurer-Cartan form of $\, \x$. From \eqref{4delPsi}, the
condition for the deformation to preserve $\, \chi$ is expressed by
the  single linear equation \beq\label{43delcon} v_1=0. \eeq We wish
to  give an  analysis of   the  compatibility equations for  the
$\, \Psi$-deformation $\, \del$  derived from \eqref{43delcon}.

Differentiating \eqref{43delcon}, one gets
\beq\label{43delstructure}
w_1= u_{{0}}
(-3\,u_{{0}}+2\,a_{{1}}+2\,a_{{2}}+2\,a_{{4}}).
\eeq
Since $\, v_1=0$ and $\, w_1$ is a function of $\, u_0$,
\emph{there exists at most one parameter family of
$\,(\Psi, \, \chi)$-deformations}.
Differentiating the equation for $\,  w_1$ again, one gets the
integrability equation
\beq
u_0 \, (a_{{4 1}}-4\,b_{{1}}+2\,b_{{2}}) =0 . \n
\eeq
If $\, a_{{4 1}}-4\,b_{{1}}+2\,b_{{2}}\ne 0$, this
forces $\, u_0=0$  and the deformation is trivial. Hence we must
have
\beq\label{43uni}
a_{{4 1}}-4\,b_{{1}}+2\,b_{{2}}=0.
\eeq
\begin{remark}
A similar analysis shows that for a general nondegenerate Legendrian
surface, either it admits maximum one parameter family of
$(\Psi, \, \chi)$-deformations, or  it does not admit any such deformations.
The Legendrian surfaces which admit  $(\Psi, \, \chi)$-deformations
are characterized by the following two relations on the structure
coefficients;
\begin{align}
a_2-a_3&= 0, \n \\
a_{{4 1}}-a_{12} &=4\,b_{{1}}. \n
\end{align}
\end{remark}

Successively differentiating \eqref{43uni}, one gets a set of three
compatibility equations.
\begin{align}\label{43key}
b_1&=0,   \quad b_{11}=0,\\
( a_{{1}}-a_{{4}}) b_{{2}}&=0,  \n  \\
( a_{{1}}-a_{{4}}) (c_1+a_1\, a_4-a_2^2)&=0. \n
\end{align}
\noi At this juncture, the analysis divides into two cases.

\texttt{Case  $\, a_1-a_4=0$}. This is a subset   of
$D_0$-surfaces. It is easily checked that the structure equation
\eqref{R0st} remains in involution with the additional condition
$\, b_1=b_{11}=0$.

\begin{definition}\label{S0defi}
A \emph{Legendrian $S_0$-surface} is a  Legendrian $D_0$-surface for
which the structure coefficients satisfy the additional relation
$\, b_1=b_{11}=0$.
\end{definition}

Let us record the full structure equation  for $S_0$-surfaces.
\begin{align}\label{S0st}
da_1+\frac{4}{3} \,a_{{1}}\omega_{{0 0}}&=
-2\,b_{{2}} ( \omega^{{1}}+\omega^{{2}}),   \\
da_2+\frac{4}{3}\,a_{{2}}\omega_{{0 0}}&=
3\,b_{{2}} ( \omega^{{1}}+\omega^{{2}}),    \n \\
db_2+2\,b_{{2}}\omega_{{0 0}}&=
 - (c_{{1}}+a_{{1}}^2-{a_{{2}}}^{2})( \omega^{{1}}+\omega^2), \n \\
dc_{1}+\frac{8}{3}\,c_{{1}}\omega_{00}&= c_{{1 1}}
(\omega^{{1}}+\omega^2).  \n
\end{align}
\noi The induced Maurer-Cartan form $\, \phi$ takes the following
form. \beq \phi =
\begin{pmatrix}
\omega_{{0 0}}& (a_1+a_2) \omega^1&(a_1+a_2) \omega^2
&c_{{1}}(\omega^1+\omega^2) & b_{{2}}\omega^2 & b_{{2}}  \omega^1
\\\noalign{\medskip}\omega^1&\frac{1}{3} \,\omega_{{0 0}}&\cdot
& b_{{2}}\omega^2&a_{{1}}\omega^2 &a_{{2}}(\omega^1+\omega^2)
\\\noalign{\medskip}\omega^2&\cdot&\frac{1}{3}\,\omega_{{0 0}}
& b_{{2}} \omega^1 &a_{{2}}(\omega^1+\omega^2)&a_{{1}}\omega^1
\\\noalign{\medskip}\cdot&\cdot&\cdot&-\omega_{{0 0}}&-\omega^1&-\omega^2
\\\noalign{\medskip}\cdot&\omega^2&\omega^1+\omega^2
&-(a_1+a_2) \omega^1&-\frac{1}{3} \,\omega_{{0 0}}&\cdot
\\\noalign{\medskip}\cdot&\omega^1+\omega^2&\omega^1
&-(a_{{1}}+a_2) \omega^2&\cdot&-\frac{1}{3} \,\omega_{{0 0}}
\end{pmatrix}.
\eeq Note that the subset of $\, \Psi$-null surfaces with the
structure coefficient $\, b=0$, and   tri-ruled surfaces are
examples of  Legendrian $S_0$-surfaces.

\texttt{Case  $\, a_1-a_4 \ne 0$}. This is a subset of
$D$-surfaces. It is easily checked that the structure equation
\eqref{R0st} remains compatible with the additional condition $\,
b_1=b_2=a_{41}=b_{11}=0$, $\, c_1=a_2^2-a_1\,a_4$.

\begin{definition}\label{Sdefi}
A \emph{Legendrian $S$-surface} is a  Legendrian $D$-surface for
which the structure coefficients satisfy the additional relation $\,
b_1=b_2=a_{41}=b_{11}=0$, and $\, c_1=a_2^2-a_1 a_4$.
\end{definition}

Let us record the full structure equation  for $D$-surfaces.
\begin{align}\label{Sst}
da_1+\frac{4}{3} \,a_{{1}}\omega_{{0 0}}&=0,   \\
da_2+\frac{4}{3}\,a_{{2}}\omega_{{0 0}}&=0,  \n \\
da_4+\frac{4}{3}\,a_{{4}}\omega_{{0 0}}&=0. \n
\end{align}
The induced Maurer-Cartan form $\, \phi$ takes the following form.
\beq \phi =
\begin{pmatrix}
\omega_{{0 0}}& (a_4+a_2) \omega^1&(a_1+a_2) \omega^2
&(a_2^2-a_1\,a_4)(\omega^1+\omega^2) &\cdot &\cdot
\\\noalign{\medskip}\omega^1&\frac{1}{3} \,\omega_{{0 0}}&\cdot
&\cdot&a_{{1}}\omega^2 &a_{{2}}(\omega^1+\omega^2)
\\\noalign{\medskip}\omega^2&\cdot&\frac{1}{3}\,\omega_{{0 0}}
&\cdot &a_{{2}}(\omega^1+\omega^2)&a_{{4}}\omega^1
\\\noalign{\medskip}\cdot&\cdot&\cdot&-\omega_{{0 0}}&-\omega^1&-\omega^2
\\\noalign{\medskip}\cdot&\omega^2&\omega^1+\omega^2
&-(a_4+a_2) \omega^1&-\frac{1}{3} \,\omega_{{0 0}}&\cdot
\\\noalign{\medskip}\cdot&\omega^1+\omega^2&\omega^1
&-(a_{{1}}+a_2) \omega^2&\cdot&-\frac{1}{3} \,\omega_{{0 0}}
\end{pmatrix}.
\eeq Note that when $\, a_4=a_1$, the structure equation for $\,
S$-surfaces degenerates to  the structure equation for $\,
S_0$-surfaces with the additional condition $\, b_2=0, \,
c_1=-a_1^2+a_2^2$.

\begin{theorem}\label{43thm}
The set of maximally $\, \Psi$-deformable Legendrian surfaces in $\,
\PP^5$ which admit one parameter family of $\,( \Psi, \,
\chi)$-deformations fall  into two categories; Legendrian $\,
S_0$-surfaces, or Legendrian $\, S$-surfaces. A general Legendrian
$\, S_0$-surface depends on one arbitrary function of 1 variable,
whereas a  general Legendrian $\, S$-surface depends on  two
constants.
\end{theorem}
\Pf The  structure equation for $\, S_0$-surfaces is in involution
with the last nonzero Cartan character s$_1=1$(we omit the details).
For  $\, S$-surfaces, consider  the invariant map
$\, I=(a_1, \, a_2, \, a_4):$ F$\to \C^3$,
where F   is the canonical bundle of 5-adapted
frames. Since $\, I$ generically has rank 1, the local moduli space
of $\, S$-surfaces has general dimension
dim $(\C^3) -$ rank$(I)=2$. $\sq$

\section{Examples}\label{sec5}
In this final section, we give a differential geometric
characterization of  tri-ruled surfaces,  Section \ref{sec241},
which are examples of  Legendrian $S_0$-surfaces. In Section
\ref{sec51}, the flat case is characterized as a part of a
Legendrian map from $\, \PP^2$ blown up at three distinct collinear
points. In Section \ref{sec52}, the non-flat case is characterized
as a part of a Legendrian embedding from $\, \PP^2$ blown up at
three non-collinear points. In both cases, the Legendrian map is
given by a system of cubics through the three points.

Let $\, (X_0, X_1,  X_2, Y_0,  Y_1,  Y_2)$ be  the standard adapted
coordinate of $\, \C^6$ such that the symplectic 2-form $\, \varpi =
dX_0 \w dY_0 + dX_1 \w dY_1 +dX_2 \w dY_2$.
\subsection{Flat surface}\label{sec51}
This is the class of surface for which all   the structure
coefficients vanish; $\, a_i=b_j=c_k=0$.

Since $\, d\omega_{00}=0$, take a section of the frame for which $\,
\omega_{00}=0$, and  $\, d\omega^1=d\omega^2=0$ consequently.
Introduce a local coordinate $\, (x, \, y)$ such that $\,
\omega^1=dx, \, \omega^2=dy$, and  express the Maurer-Cartan form
$\, \phi = A\, dx + B\,dy$ for constant coefficient matrices $\, A,
\, B$. Since $\, d\phi= - \phi \w \phi =0$, $\, A$ and $\, B$
commute and $\, g = \exp(A x+B y)$ is a solution of  the defining
equation \beq g^{-1} \, dg = \phi. \n \eeq The exponential can be
computed, and by definition of $\, \phi$ in Section \ref{sec2}, the
first column of $\, g$ gives the following local parametrization of
the flat Legendrian surface. \beq \x_1(x, y)=
\begin{pmatrix}
X_0 \\ X_1 \\ X_2 \\Y_0 \\ Y_1 \\ Y_2
\end{pmatrix}
=
\begin{pmatrix}
1 \\ x \\y \\ -\frac{x y}{2}(x+y) \\  x y + \frac{y^2}{2} \\ x y +
\frac{x^2}{2}
\end{pmatrix}. \n
\eeq

\begin{theorem}\label{thm51}
Let  $\, \pi: M_3  \to \PP^2$ be the rational surface obtained by
blowing up  $\, \PP^2$  at three distinct collinear points $\, \{ \,
p_1, \, p_2, \, p_3 \, \}$. Let $\, L$ be the line through   $\,
p_i$'s, $\, E_i = \pi^{-1}(p_i)$ be the exceptional divisor, and let
$\, H$ be the linear divisor of $\, \PP^2$. Let $\, \widehat{L}$ be
the -2-curve, the proper transform of  $\, L$. There exists a six
dimensional proper subspace $\, W$ of the linear system $\,
[\pi^*(3H)-E_1-E_2-E_3]$ which gives a Legendrian map  $\,
\widehat{\x}: M_3 \to \PP^5$.
 $\, \widehat{\x}$  is an embedding on $\, M_3 - \widehat{L}$, and
it degenerates to a point on $\, \widehat{L}$.

\textnormal{a)} A flat Legendrian surface is locally equivalent to a
part of $\; \widehat{\x}(M_3- \cup E_i )$.

\textnormal{b)} The system $\, W$ for $\,  \widehat{\x}$ is a six
dimensional subspace of the proper transform of the set of cubics
through   $\, p_i$'s. Each -1-curve $\, E_i$ is mapped to a line
under $\, \widehat{\x}$.

\textnormal{c)} The asymptotic web is given by the  proper transform
of the three pencils of lines through $\, p_i$'s.

\end{theorem}

Proof of theorem is presented below in four steps.

\noi \texttt{Step 1.} Consider the birational map $\, \x:  \PP^2
\dashrightarrow   \PP^5$ associated  to $\, \x_1(x, y)$ defined by
\beq\label{birational51} \x([x, y, z])=
\begin{pmatrix}
X_0 \\ X_1 \\ X_2 \\Y_0 \\ Y_1 \\ Y_2
\end{pmatrix}
=
\begin{pmatrix}
z^3 \\ x z^2 \\y  z^2 \\ -\frac{x y}{2}(x+y) \\  (x y +
\frac{y^2}{2}) z \\ ( x y +   \frac{x^2}{2})z
\end{pmatrix},
\eeq where $\, [x,y,z]$ is the standard projective coordinate of $\,
\PP^2$. It is undefined at three points
\begin{align}
p_1&=[1,0,0],  \n \\
p_2&=[0,1,0],  \n \\
p_3&=[1,-1,0].  \n
\end{align}

At $\, p_1$, introduce the parametrization of the blow up by $\, [x,
\,y, \,z]=[1, \, \lambda_1   z, \, z]$ for the blow up parameter $\,
\lambda_1$. The birational map becomes \beq\label{blowup1} \x([x, y,
z])=
\begin{pmatrix}
X_0 \\ X_1 \\ X_2 \\Y_0 \\ Y_1 \\ Y_2
\end{pmatrix}
=
\begin{pmatrix}
z^3 \\  z^2 \\   \lambda_1     z^3 \\ -\frac{ \lambda_1   z}{2}(1+ \lambda_1   z) \\
( \lambda_1   + \frac{ \lambda_1^2   z}{2}) z^2 \\ (  \lambda_1   z
+   \frac{1}{2}) z
\end{pmatrix}
=
\begin{pmatrix}
z^2 \\  z  \\   \lambda_1     z^2 \\ -\frac{ \lambda_1 }{2}(1+ \lambda_1   z) \\
( \lambda_1   + \frac{ \lambda_1^2   z}{2}) z \\ (  \lambda_1   z +
\frac{1}{2})
\end{pmatrix}
\to
\begin{pmatrix}
0 \\  0  \\  0 \\ -\frac{ \lambda_1 }{2}  \\
0 \\ \frac{1}{2}
\end{pmatrix}, \quad \mbox{as}\;\;  z \to 0.
\eeq Similar formulae for $\, p_2, \, p_3$ show that the exceptional
divisors $\, E_2, \, E_3$ are mapped to \beq\label{blowup23} E_2
\to
\begin{pmatrix}
0 \\  0  \\  0 \\ -\frac{ \lambda_2}{2}  \\
 \frac{1}{2} \\ 0
\end{pmatrix}, \quad
E_3    \to
\begin{pmatrix}
0 \\  0  \\  0 \\  \frac{ \lambda_3}{2}  \\
 -\frac{1}{2}   \\ -\frac{1}{2}
\end{pmatrix}.
\eeq

Let $\, L = \{ \, z=0 \, \} \subset \PP^2$ be the line through $\,
p_i$'s. By definition,  $\, \x(L) =   [0,0,0,1,0,0] = \x_0  \in
\PP^5$. One may check  that $\, \x: \PP^2-L     \to \PP^5$ is an
embedding, and that the image   $\, \x( \PP^2 - L)$ is disjoint from
the exceptional loci \eqref{blowup1}, \eqref{blowup23}.

A computation with \eqref{blowup1} at $\, p_1$,  and similar
computations   at $\, p_2, \, p_3$ show  that the associated  lift
$\, \widehat{\x}: M_3 \to \PP^5$ is well defined and holomorphic,
and that $\, \widehat{\x}: M_3 - \widehat{L}  \to \PP^5$ is a smooth
embedding.

\noi \texttt{Step 2.} Consider alternatively the following
polynomial equations satisfied by $\, \x$.
\begin{align}
3\,X_0 Y_0 + X_1 Y_2+ X_2 Y_2 &=0, \n \\
X_0^2 Y_0 + \frac{1}{2} \, X_1  X_2 (X_1+X_2) &= 0,  \n \\
X_0 Y_1 -  X_1  X_2 - \frac{1}{2} X_2^2 &=0, \n \\
X_0 Y_2 -  X_1  X_2 - \frac{1}{2} X_1^2 &=0,\n  \\
\left(  X_1(Y_2-\frac{1}{2}Y_1) - X_2(Y_1-\frac{1}{2}Y_2) \right)
Y_0 &= Y_1 Y_2(Y_1-Y_2). \n
\end{align}
By a direct computation, one can verify that this set of  equations
have rank 3  on $\,  \widehat{\x}(M_3)=\overline{\x( \PP^2 - L)}$,
except at $\, \x_0$. One can  also check that at $\, \x_0$, $\,
\widehat{\x}(M_3)$ is not smooth and  has a second order branch type
singularity(we omit the details. Note that $\, \widehat{L}$ is a
-2-curve and it cannot be blown down).

\noi \texttt{Step 3.} Let $\, D$  the hyperplane section $\, D =
\widehat{\x}^{-1}(\{ \, Y_0=0 \, \})$. From  \eqref{birational51},
the divisor consists of the proper transform of    three lines $\,
\{ \, L_1, \, L_2, \, L_3 \, \} = \{ \, y=0,\, \mbox{or} \; x=0,\,
\mbox{or} \; x+y=0 \, \} \subset \PP^2$. Hence the linear system
\beq [ D ] =  [\sum_i( \pi^*(H)-E_i)] = [\pi^*(3H)- E_1-E_2-E_3]. \n
\eeq
By definition,
$\, W$ is a subspace of the linear system $\,
[\pi^*(3H) - E_1-E_2-E_3]$, the proper transform of  cubics through
$\, p_i$'s.
Finally, $\, \langle D, \, E_i \rangle = 1$, and
each $\,E_i$ is mapped to a line.

\noi \texttt{Step 4.}
One may check by direct computation that
the asymptotic web is given by the foliations
$\,  dy=0, \, dx=0, \, dx+dy=0 $, which represent  three pencils of
lines through $\, p_1, \, p_2, \, p_3$ respectively.

\subsubsection{Generalization}\label{sec511}
The construction of  flat surface admits a straightforward
generalization.

Let $\, f_k(x,y)$ be a homogeneous polynomial of degree $\, k$ for
$\, k = 3, \, 4, \, ... \, m$, such that the top degree $\, f_m(x,y)$
has no multiple factors (product of $\, m$ mutually non-proportional
linear functions in $\, x, \, y$). Consider the associated
birational map $\, \x:  \PP^2  \dashrightarrow   \PP^5$ defined by
\beq \x([x, y, z])=
\begin{pmatrix}
X_0 \\ X_1 \\ X_2 \\Y_0 \\ Y_1 \\ Y_2
\end{pmatrix}
=
\begin{pmatrix}
z^m \\ x z^{m-1} \\y  z^{m-1} \\ - \sum_{k=3}^m \, f_k \, z^{m-k}
\\ \sum_{k=3}^m \, \frac{1}{k-2} (\frac{\partial}{\partial x} f_k) \, z^{m-k+1}
\\ \sum_{k=3}^m \, \frac{1}{k-2} (\frac{\partial}{\partial y} f_k) \, z^{m-k+1}
\end{pmatrix}. \n
\eeq A direct computation shows that $\, \x$ is Legendrian.

$\, \x$ is undefined at $\, m$ points $\, \{ \, z=0, \, f_m(x,y)=0
\}  \subset \PP^2$. Let  $\, \pi: M_m  \to \PP^2$ be the rational
surface obtained by blowing up  $\, \PP^2$  at these points.
Let $\,E_i=\pi^{-1}(p_i), \, i = 1, \, 2, \, .. m$,
be the exceptional divisor, and
let $\,\widehat{L}$ be the proper transform of
the line $\, \{ \, z=0 \,\}$.
An analysis similar as above shows that $\, \x$ admits
a well defined smooth  lift
$\, \widehat{\x}: M_m \to \PP^5$.
But when $\, m \geq 4 $,
the singular locus of  $\, \widehat{\x}$ consists of $\, \widehat{L}$,
and one point from each $\, E_i - \widehat{L}$.

This class of Legendrian surfaces were first introduced  in
\cite{Bu1}.

\subsection{Tri-ruled surface}\label{sec52}
This is the class of surface with the induced Maurer-Cartan form
\eqref{23triphi}. Since the flat case is already treated,
we examine the case  $\, a \ne 0$.

From \eqref{24da}, one may scale $\, a=1$. Then $\, \omega_{00}=0$,
and  $\, d\omega^1=d\omega^2=0$ consequently. Introduce a local
coordinate $\, (s, \, t)$ such that $\,
\omega^1=\frac{\im}{\sqrt{2}} ds, \, \omega^2=\frac{\im}{\sqrt{2}}
dt$, where $\, \im^2=-1$, and  express the Maurer-Cartan form $\,
\phi = A\, ds + B\,dt$ for constant coefficient matrices $\, A, \,
B$. As in Section \ref{sec51}, the equation $\, g^{-1} dg = \phi$
can be integrated and  one gets the following local parametrization
of a  tri-ruled Legendrian surface up to conformal symplectic
transformation. \beq \x_1(s, t)=
\begin{pmatrix}
X_0 \\ X_1 \\ X_2 \\Y_0 \\ Y_1 \\ Y_2
\end{pmatrix}
=
\begin{pmatrix}
\sin(s+t)  \\ \sin(s) \\ \sin(t)  \\
 -\cos(s+t)  \\ \cos(s) \\ \cos(t)
\end{pmatrix}. \n
\eeq By a conformal symplectic transformation, we mean a linear
transformation of  $\, \C^6$ that preserves the symplectic form up
to nonzero scale, e.g., a linear transformation $\, (X, \, Y) \to (
l_1 X, \,  l_2 Y)$  for nonzero $\,l_1, \,l_2$.

An analysis shows that this local parametrization gives rise to a
Legendrian embedding of $\, \PP^1\times \PP^1 = Q^1 \times Q^1
\subset \PP^2 \times \PP^2$, the product of two conics, blown up at
two points. Since this surface is isomorphic to $\, \PP^2$ blown up
at three non-collinear points, consider the associated Legendrian
birational map $\, \x:  \PP^2  \dashrightarrow   \PP^5$ defined by
\beq\label{birational52} \x([x, y,z])=
\begin{pmatrix}
X_0 \\ X_1 \\ X_2 \\Y_0 \\ Y_1 \\ Y_2
\end{pmatrix}
=
\begin{pmatrix}
(x^2-y^2) z \\  (y^2-z^2) x  \\ (z^2-x^2) y  \\
(x^2+y^2) z \\  (y^2+z^2) x  \\  (z^2+x^2) y
\end{pmatrix},
\eeq where $\, [x,y,z]$ is the standard projective coordinate of $\,
\PP^2$,  Lemma \ref{521lemm}. It is undefined at three points
\begin{align}\label{52three}
p_1&=[1,0,0],    \\
p_2&=[0,1,0],  \n \\
p_3&=[0,0,1].  \n
\end{align}
\begin{theorem}\label{thm52}
Let  $\, \pi:  N_3  \to \PP^2$ be the rational surface obtained by
blowing up  $\, \PP^2$  at three  non-collinear points
$\, \{ \,p_1, \, p_2, \, p_3 \, \}$. Let $\, L_{k}$ be the line through
$\,(p_i, \, p_j)$, $\,(ijk)=(123)$, $\, E_i = \pi^{-1}(p_i)$ be the
exceptional divisor, and let $\, H$ be the linear divisor of
$\,\PP^2$. Let $\, \widehat{L}_{k}$ be the proper transform of
$\,L_{k}$. There exists a six dimensional proper subspace $\, W$ of the
linear system $\, [\pi^*(3H)-E_1-E_2-E_3]$ which gives a Legendrian
embedding  $\, \widehat{\x}: N_3 \hook \PP^5$.

\textnormal{a)} A non-flat tri-ruled Legendrian surface is locally
equivalent to a part of
$\; \widehat{\x}(N_3- \cup E_i \cup \widehat{L}_i)$.

\textnormal{b)} The system $\, W$ for $\, \widehat{\x}$ is a six
dimensional subspace of the proper transform of the set of cubics
through   $\, p_i$'s. Each of the six -1-curves $\, E_i$  and
$\,\widehat{L}_{k}$ is mapped to a line under the embedding.

\textnormal{c)} The asymptotic web is given by the  proper transform
of the three pencils of lines through $\, p_i$'s.

\end{theorem}

Proof of theorem is presented below in four steps.

\noi \texttt{Step 1.} $\, \widehat{\x}$ is an immersion:

By a direct computation, it is verified that $\, \x$ is an immersion
on  $\, \PP^2 - \cup p_i$.

At $\, p_1$, introduce the parametrization of the blow up by $\, [x,
\,y, \,z]=[1, \, \lambda_1   z, \, z]$ for the blow up parameter $\,
\lambda_1$. The birational map becomes \beq\label{blowup521} \x([x,
y, z])=
\begin{pmatrix}
X_0 \\ X_1 \\ X_2 \\Y_0 \\ Y_1 \\ Y_2
\end{pmatrix}
=
\begin{pmatrix}
(x^2-y^2) z \\  (y^2-z^2) x  \\ (z^2-x^2) y  \\
(x^2+y^2) z \\  (y^2+z^2) x  \\  (z^2+x^2) y
\end{pmatrix}
=
\begin{pmatrix}
(1-\lambda_1^2 z ^2)   \\  (\lambda_1^2  -1)   z  \\ (z^2-1) \lambda_1    \\
(1+\lambda_1^2 z ^2)   \\  (\lambda_1^2  +1)   z\\  (z^2+1)
\lambda_1
\end{pmatrix}
\to
\begin{pmatrix}
1 \\  0  \\  - \lambda_1\\
1 \\  0  \\     \lambda_1
\end{pmatrix}, \quad \mbox{as}\;\;  z \to 0.
\eeq Similar computations for $\, p_2, \, p_3$ show that $\, \x$
admits a lift  $\, \widehat{\x}: N_3 \to \PP^5$ such that $\,
\widehat{\x}(N_3) = \overline{\x(\PP^2 - \cup p_i)}$. The
exceptional divisors are respectively mapped to
\beq\label{blowupp23} E_2    \to
\begin{pmatrix}
-\lambda_2  \\  1  \\  0 \\ \lambda_2
\\ 1 \\ 0
\end{pmatrix}, \quad
E_3    \to
\begin{pmatrix}
0 \\  -\lambda_3  \\  1 \\
0 \\   \lambda_3  \\  1
\end{pmatrix}, \n
\eeq for blow up parameters $\, \lambda_2, \, \lambda_3$.

From \eqref{blowup521}, one may  check that the three vectors $\,
\x, \, \frac{\partial \x}{\partial z}, \, \frac{\partial
\x}{\partial \lambda_1}$ are independent at $\, z=0$. Similar
computations for $\, p_2, \, p_3$ show that $\,  \widehat{\x}$ is an
immersion on the exceptional divisors $\, E_i$. Hence  $\,
\widehat{\x}$ is an immersion on $\, N_3$.

\noi \texttt{Step 2.}  $\, \widehat{\x}$ is injective:

It is clear that  $\, \widehat{\x}$ is injective on $\, \cup E_i$,
and that $\, \widehat{\x}(\cup E_i)$ is disjoint from $\, \x(\PP^2-
\cup p_i)$. It suffices to show that $\, \x$ is injective on $\,
\PP^2- \cup p_i$.

Suppose $\, \x([x, y, z]) = \x([x', y', z'])$. Then
\begin{align}
x^2y&=\mu \,x'^2y', \quad x^2z =\mu \,x'^2z', \n \\
y^2z&=\mu \,y'^2z', \quad y^2x =\mu \,y'^2x',  \n  \\
z^2x&=\mu \,z'^2x', \quad z^2y =\mu \,z'^2y', \n
\end{align}
for a nonzero scaling parameter $\, \mu$.

Case  $\, x=0;  \, y, \, z \ne 0$. Then $\, x'y' = x'z'=0$. If $\,
x' \ne 0$, then $\, y'z'=0 = yz$, a contradiction. Hence $\, x'=0$.
The remaining equations then show that $\, \frac{y}{z} =
\frac{y'}{z'}$.

Case  $\, x,\, y,\, z \ne 0$.  One has $\, \frac{x}{y} =
\frac{x'}{y'}, \, \frac{y}{z} =  \frac{y'}{z'}, \, \frac{z}{x} =
\frac{z'}{x'}$. Hence $\, [x, y, z] = [x', y', z']$.

\noi \texttt{Step 3.} Let $\, D$  the hyperplane section $\, D =
\widehat{\x}^{-1}( \{ \, X_0+X_1+X_2=0 \, \})$. From
\eqref{birational52}, $\, X_0+X_1+X_2=(x-y)(y-z)(z-x)$. Hence the
linear system \beq [ D ] =  [\sum_i( \pi^*(H)-E_i)] = [\pi^*(3H)-
E_1-E_2-E_3]. \n \eeq $\, W$ is a subspace of the linear system of
the proper transform of  cubics through $\, p_i$'s,
and $\, \widehat{\x}$ is not normal.

Since $\, [\widehat{L}_{k}] = [\pi^*(H) - E_i -E_j], \,
(ijk)=(123)$, one has $\, \langle D, \, E_i \rangle =  \langle D,
\,\widehat{L}_{k} \rangle = 1$, and  each  $\, E_i$  and  $\,
\widehat{L}_{k}$ is mapped to a line.

\noi \texttt{Step 4.} The equations for the asymptotic web can be
checked on the affine chart $\, [x, y, 1]$ by a direct computation.
We omit the details. Let $\, C$ be a line on $\, \PP^2$ that passes
through exactly one of $\, p_i$'s. The proper transform $\,
\widehat{C}$ of $\, C$ has the divisor class $\, \pi^*(H)-E_i$.
Hence $\, \langle D, \, \widehat{C} \rangle = 2$, and  each leaf of
the asymptotic foliations  is mapped to a linear
$\, \PP^2 \subset \PP^5$ which is necessarily Legendrian
from the defining properties of the tri-ruled surfaces.
$\sq$

\begin{remark}
The linear system of  conics through three non-collinear points
gives the classical quadratic transformation of $\, \PP^2$.
\end{remark}

As the three points degenerate to become collinear, $\,
\widehat{\x}(N_3)$ degenerates to the flat surface in Section
\ref{sec51}. The isolated singularity of the flat Legendrian surface
thus admits a smoothing.

It is not known  if  every  del Pezzo surface admits a Legendrian
embedding. Legendrian embeddings of a set of degree 4 del Pezzo
surfaces were constructed in \cite{Bu2}.

Note the algebraic equations satisfied by $\,  \widehat{\x}$(this is
not a complete intersection).
\begin{align}\label{52cone}
&X_0^2-Y_0^2 =X_1^2-Y_1^2=X_2^2-Y_2^2,   \\
& X_0 (X_1X_2+Y_1Y_2)+Y_0( X_1Y_2+Y_1X_2)  =0. \n
\end{align}
Legendrian surface $\, \widehat{\x}(N_3)$ can thus be considered as
a complexification  of the homogeneous special Legendrian torus with
parallel second fundamental form, \cite{HL}.


\subsubsection{Generalization}\label{sec521}
The construction of  tri-ruled surface admits a straightforward
generalization.

Let $\, m, \, n$ be a pair of  positive integers.
Let $\, f_{m}, \, g_{m}$ be the homogeneous polynomials of
degree $\, m$ of  two variables which represent  $\, \sin(m s ), \,
\cos(m s)$;
\begin{align}
\sin(m s)&=f_{m }( \sin(s),    \, \cos(s)) =
\sum_{0 \leq k \leq \frac{m-1}{2}}  (-1)^k {m \choose 2k+1} \sin^{2k+1}(s) \cos^{m-(2k+1)}(s) ,  \n \\
\cos(m s )&=g_{m }( \sin(s), \,    \cos(s)) = \sum_{0 \leq k \leq
\frac{m}{2}}  (-1)^k {m \choose 2k} \sin^{2k}(s) \cos^{m-2k}(s) . \n
\end{align}
Consider the following local parametrization of  a Legendrian
surface. \beq\label{521local} \x_1(s, t)=
\begin{pmatrix}
X_0 \\ X_1 \\ X_2 \\Y_0 \\ Y_1 \\ Y_2
\end{pmatrix}
=
\begin{pmatrix}
 \sin(m s + n t)  \\ \sqrt{m} \sin(s)    \\ \sqrt{n}\sin(t)  \\
 - \cos(m s + n t)  \\ \sqrt{m}\cos(s) \\ \sqrt{n}\cos(t)
\end{pmatrix}
=
\begin{pmatrix}
f_m(u_1, \, u_2) g_n(v_1, \, v_2) + g_m(u_1, \, u_2) f_n(v_1, \, v_2) \\
(-\im)^{m+n-1} \sqrt{m} \, u_1 u_0^{m-1} v_0^n  \\
(-\im)^{m+n-1}\sqrt{n} \, v_1 u_0^m   v_0^{n-1}  \\
f_m(u_1, \, u_2) f_n(v_1, \, v_2) -  g_m(u_1, \, u_2) g_n(v_1, \, v_2) \\
(-\im)^{m+n-1}\sqrt{m}\, u_2 u_0^{m-1} v_0^n \\
(-\im)^{m+n-1}\sqrt{n} \, v_2 u_0^m   v_0^{n-1}
\end{pmatrix},
\eeq where $\,  u_0^2+u_1^2+u_2^2=0, \, v_0^2+v_1^2+v_2^2=0$.

\begin{lemma}\label{521lemm}
The local parametrization \eqref{521local} is equivalent to the
following Legendrian birational map $\, \x:  \PP^2  \dashrightarrow
\PP^5$ up to conformal symplectic transformation.
\beq\label{birational521}
\x([x, y, z])=
\begin{pmatrix}
X_0 \\ X_1 \\ X_2 \\Y_0 \\ Y_1 \\ Y_2
\end{pmatrix}
=
\begin{pmatrix}
x^{2m} z  - y^{2n} z^{2m-2n+1}  \\
\sqrt{m} (z^2-x^2) x^{m-1} y^n   z^{m-n} \\
\sqrt{n}  (y^2-z^2) x^{m}    y^{n-1}   z^{m-n} \\
x^{2m} z   +  y^{2n} z^{2m-2n+1}  \\
\sqrt{m} (z^2+x^2) x^{m-1} y^n   z^{m-n} \\
\sqrt{n}  (y^2+z^2) x^{m}    y^{n-1}   z^{m-n}
\end{pmatrix}.
\eeq
\end{lemma}
\Pf Let $\, [x, y,z]$ be the homogeneous coordinate of $\, \PP^2$.
Take the following birational map
$\, \varphi: \PP^2 \dashrightarrow
Q^1\times Q^1  \subset \PP^2 \times \PP^2$;
\begin{align}
\varphi([x,y,z]) &= ( [u_0, \, u_1, \, u_2], \, [v_0, \, v_1, \, v_2] ), \n \\
                         &= ( [2xz, \, z^2-x^2, \, \im (z^2+x^2)], \, [2yz, \, -z^2+y^2, \, \im (z^2+y^2)]).\n
\end{align}
Lemma follows from  de Moivre's formula, \beq (u_2 \pm \im u_1)^m
(v_2 \pm \im v_1)^n = (g_m(u_1, u_2) \pm  \im f_m(u_1, u_2))(
g_n(v_1, v_2) \pm \im f_n(v_1, v_2)). \n  \sq \eeq

Consider the case $\, m=n$.
\eqref{birational521} is undefined at the three points  of \eqref{52three}.
Let $\,\widehat{ \x }: N_3 \hook \PP^5$ be the induced lift.
It is never an immersion  except when $\, (m,n)  =  (1,1)$.
For example when $\,m \geq 2$,
the line $\, \{ \, y=0 \, \}$ degenerates to a point.



 \vsp{1pc}

\end{document}